% TO GET TABLE OF CONTENT, TEX TWICE.

%%%%%%%%%%%%%%%%%%  tex macros for preprints, cm version %%%%%%%%%%%%%%
%                     (P. Ginsparg, last updated 9/91)
%                if confused, type `b' in response to query 
%
%---------------------------------------------------------------------%
%% site dependent options: 
%% \unredoffs and \redoffs define horizontal and vertical offsets 
%% respectively for unreduced and reduced modes. \speclscape defines
%% the \special{} call that sets printer to landscape (sideways) mode.
%% from standard set below, leave uncommented as appropriate or redefine
%
%%% next 400dpi
%\def\unredoffs{} \def\redoffs{\voffset=-.31truein\hoffset=-.48truein}
%\def\speclscape{\special{landscape}}
%
%%% apple lw
\def\unredoffs{} \def\redoffs{\voffset=-.31truein\hoffset=-.59truein}
\def\speclscape{\special{ps: landscape}}
%
%%% qms lasergrafix:
%\def\unredoffs{} \def\redoffs{\voffset=-.4truein\hoffset=.125truein}
%\def\speclscape{\special{qms: landscape}}
%
%%% saclay A4 paper:
%\def\unredoffs{\hoffset-.14truein\voffset-.2truein} 
%\def\redoffs{\voffset=-.45truein\hoffset=-.21truein} 
%\def\speclscape{\special{landscape}}
%
%---------------------------------------------------------------------%
%
\newbox\leftpage \newdimen\fullhsize \newdimen\hstitle \newdimen\hsbody
\tolerance=1000\hfuzz=2pt
\catcode`\@=11 % This allows us to modify PLAIN macros.
%\def\bigans{b }    %comment two lines to remove option
%\message{ big or little (b/l)? }\read-1 to\answ
%
\ifx\answ\bigans\message{(This will come out unreduced.}
\magnification=1200\unredoffs\baselineskip=16pt plus 2pt minus 1pt
\hsbody=\hsize \hstitle=\hsize %take default values for unreduced format
\else\message{(This will be reduced.} \let\l@r=L
\magnification=1000\baselineskip=16pt plus 2pt minus 1pt \vsize=7truein
\redoffs \hstitle=8truein\hsbody=4.75truein\fullhsize=10truein\hsize=\hsbody
\output={\ifnum\pageno=0 %%% This is the HUTP version
  \shipout\vbox{\speclscape{\hsize\fullhsize\makeheadline}
    \hbox to \fullhsize{\hfill\pagebody\hfill}}\advancepageno
  \else
  \almostshipout{\leftline{\vbox{\pagebody\makefootline}}}\advancepageno 
  \fi}
\def\almostshipout#1{\if L\l@r \count1=1 \message{[\the\count0.\the\count1]}
      \global\setbox\leftpage=#1 \global\let\l@r=R
 \else \count1=2
  \shipout\vbox{\speclscape{\hsize\fullhsize\makeheadline}
      \hbox to\fullhsize{\box\leftpage\hfil#1}}  \global\let\l@r=L\fi}
\fi
%---------------------------------------------------------------------
%
\newcount\yearltd\yearltd=\year\advance\yearltd by -1900

\def\Title#1#2{\nopagenumbers\abstractfont\hsize=\hstitle\rightline{#1}%
\vskip 1in\centerline{\titlefont #2}\abstractfont\vskip .5in\pageno=0}
\def\Date#1{\vfill\leftline{#1}\tenpoint\supereject\global\hsize=\hsbody%
\footline={\hss\tenrm\folio\hss}}% 	restores pagenumbers
%
%       use following instead of \Date on the preliminary draft, 
%       puts date/time on each page in big mode, writes labels in margins

\def\draftmode{\message{ DRAFTMODE }\def\draftdate{{\rm preliminary draft:
\number\month/\number\day/\number\yearltd\ \ \hourmin}}%
\headline={\hfil\draftdate}\writelabels\baselineskip=20pt plus 2pt minus 2pt
 {\count255=\time\divide\count255 by 60 \xdef\hourmin{\number\count255}
  \multiply\count255 by-60\advance\count255 by\time
  \xdef\hourmin{\hourmin:\ifnum\count255<10 0\fi\the\count255}}}
%       use \nolabels to get rid of eqn, ref, and fig labels in draft mode
\def\nolabels{\def\wrlabeL##1{}\def\eqlabeL##1{}\def\reflabeL##1{}}
\def\writelabels{\def\wrlabeL##1{\leavevmode\vadjust{\rlap{\smash%
{\line{{\escapechar=` \hfill\rlap{\sevenrm\hskip.03in\string##1}}}}}}}%
\def\eqlabeL##1{{\escapechar-1\rlap{\sevenrm\hskip.05in\string##1}}}%
\def\reflabeL##1{\noexpand\llap{\noexpand\sevenrm\string\string\string##1}}}
\nolabels
%
% tagged sec numbers
\global\newcount\secno \global\secno=0
\global\newcount\meqno \global\meqno=1
\def\newsec#1{\global\advance\secno by1\message{(\the\secno. #1)}
%\ifx\answ\bigans \vfill\eject \else \bigbreak\bigskip \fi  %if desired
\global\subsecno=0\eqnres@t\noindent{\bf\the\secno. #1}
\writetoca{{\secsym} {#1}}\par\nobreak\medskip\nobreak}
\def\eqnres@t{\xdef\secsym{\the\secno.}\global\meqno=1\bigbreak\bigskip}
\def\sequentialequations{\def\eqnres@t{\bigbreak}}\xdef\secsym{}
\global\newcount\subsecno \global\subsecno=0
\def\subsec#1{\global\advance\subsecno by1\message{(\secsym\the\subsecno. #1)}
\ifnum\lastpenalty>9000\else\bigbreak\fi
\noindent{\it\secsym\the\subsecno. #1}\writetoca{\string\quad 
{\secsym\the\subsecno.} {#1}}\par\nobreak\medskip\nobreak}
\def\appendix#1#2{\global\meqno=1\global\subsecno=0\xdef\secsym{\hbox{#1.}}
\bigbreak\bigskip\noindent{\bf Appendix #1. #2}\message{(#1. #2)}
\writetoca{Appendix {#1.} {#2}}\par\nobreak\medskip\nobreak}
%
%       \eqn\label{a+b=c}	gives displayed equation, numbered
%				consecutively within sections.
%     \eqnn and \eqna define labels in advance (of eqalign?)
%
\def\eqnn#1{\xdef #1{(\secsym\the\meqno)}\writedef{#1\leftbracket#1}%
\global\advance\meqno by1\wrlabeL#1}
\def\eqna#1{\xdef #1##1{\hbox{$(\secsym\the\meqno##1)$}}
\writedef{#1\numbersign1\leftbracket#1{\numbersign1}}%
\global\advance\meqno by1\wrlabeL{#1$\{\}$}}
\def\eqn#1#2{\xdef #1{(\secsym\the\meqno)}\writedef{#1\leftbracket#1}%
\global\advance\meqno by1$$#2\eqno#1\eqlabeL#1$$}
%
%			 footnotes
\newskip\footskip\footskip14pt plus 1pt minus 1pt %sets footnote baselineskip
\def\footnotefont{\ninepoint}\def\f@t#1{\footnotefont #1\@foot}
\def\f@@t{\baselineskip\footskip\bgroup\footnotefont\aftergroup\@foot\let\next}
\setbox\strutbox=\hbox{\vrule height9.5pt depth4.5pt width0pt}
\global\newcount\ftno \global\ftno=0
\def\foot{\global\advance\ftno by1\footnote{$^{\the\ftno}$}}
%
%say \footend to put footnotes at end
%will cause problems if \ref used inside \foot, instead use \nref before
\newwrite\ftfile   
\def\footend{\def\foot{\global\advance\ftno by1\chardef\wfile=\ftfile
$^{\the\ftno}$\ifnum\ftno=1\immediate\openout\ftfile=foots.tmp\fi%
\immediate\write\ftfile{\noexpand\smallskip%
\noexpand\item{f\the\ftno:\ }\pctsign}\findarg}%
\def\footatend{\vfill\eject\immediate\closeout\ftfile{\parindent=20pt
\centerline{\bf Footnotes}\nobreak\bigskip\input foots.tmp }}}
\def\footatend{}
%
%     \ref\label{text}
% generates a number, assigns it to \label, generates an entry.
% To list the refs on a separate page,  \listrefs
%
\global\newcount\refno \global\refno=1
\newwrite\rfile
%
%\def\ref{[\the\refno]\nref}  %This line is replaced by
                              %the next line when references
                              %are put alphabetically at the
                              %beginning of the paper
\def\ref{\nref}
\def\nref#1{\xdef#1{[\the\refno]}\writedef{#1\leftbracket#1}%
\ifnum\refno=1\immediate\openout\rfile=refs.tmp\fi
\global\advance\refno by1\chardef\wfile=\rfile\immediate
\write\rfile{\noexpand\item{#1\ }\reflabeL{#1\hskip.31in}\pctsign}\findarg}
%	horrible hack to sidestep tex \write limitation
\def\findarg#1#{\begingroup\obeylines\newlinechar=`\^^M\pass@rg}
{\obeylines\gdef\pass@rg#1{\writ@line\relax #1^^M\hbox{}^^M}%
\gdef\writ@line#1^^M{\expandafter\toks0\expandafter{\striprel@x #1}%
\edef\next{\the\toks0}\ifx\next\em@rk\let\next=\endgroup\else\ifx\next\empty%
\else\immediate\write\wfile{\the\toks0}\fi\let\next=\writ@line\fi\next\relax}}
\def\striprel@x#1{} \def\em@rk{\hbox{}} 
\def\lref{\begingroup\obeylines\lr@f}
\def\lr@f#1#2{\gdef#1{\ref#1{#2}}\endgroup\unskip}

\def\addref#1{\immediate\write\rfile{\noexpand\item{}#1}} %now unnecessary
\def\startrefs#1{\immediate\openout\rfile=refs.tmp\refno=#1}
\def\xref{\expandafter\xr@f}\def\xr@f[#1]{#1}
\def\refs#1{\count255=1[\r@fs #1{\hbox{}}]}
\def\r@fs#1{\ifx\und@fined#1\message{reflabel \string#1 is undefined.}%
\nref#1{need to supply reference \string#1.}\fi%
\vphantom{\hphantom{#1}}\edef\next{#1}\ifx\next\em@rk\def\next{}%
\else\ifx\next#1\ifodd\count255\relax\xref#1\count255=0\fi%
\else#1\count255=1\fi\let\next=\r@fs\fi\next}
%

%
% this is ugly, but moore insists
\newwrite\ffile\global\newcount\figno \global\figno=1
\def\fig{fig.~\the\figno\nfig}
\def\nfig#1{\xdef#1{fig.~\the\figno}%
\writedef{#1\leftbracket fig.\noexpand~\the\figno}%
\ifnum\figno=1\immediate\openout\ffile=figs.tmp\fi\chardef\wfile=\ffile%
\immediate\write\ffile{\noexpand\medskip\noexpand\item{Fig.\ \the\figno. }
\reflabeL{#1\hskip.55in}\pctsign}\global\advance\figno by1\findarg}
\def\xfig{\expandafter\xf@g}\def\xf@g fig.\penalty\@M\ {}
\def\figs#1{figs.~\f@gs #1{\hbox{}}}
\def\f@gs#1{\edef\next{#1}\ifx\next\em@rk\def\next{}\else
\ifx\next#1\xfig #1\else#1\fi\let\next=\f@gs\fi\next}
\newwrite\lfile
{\escapechar-1\xdef\pctsign{\string\%}\xdef\leftbracket{\string\{}
\xdef\rightbracket{\string\}}\xdef\numbersign{\string\#}}

\def\writestop{\def\writestoppt{\immediate\write\lfile{\string\pageno%
\the\pageno\string\startrefs\leftbracket\the\refno\rightbracket%
\string\def\string\secsym\leftbracket\secsym\rightbracket%
\string\secno\the\secno\string\meqno\the\meqno}\immediate\closeout\lfile}}
\def\writestoppt{}\def\writedef#1{}
\def\seclab#1{\xdef #1{\the\secno}\writedef{#1\leftbracket#1}\wrlabeL{#1=#1}}
\def\subseclab#1{\xdef #1{\secsym\the\subsecno}%
\writedef{#1\leftbracket#1}\wrlabeL{#1=#1}}
\newwrite\tfile \def\writetoca#1{}
\def\leaderfill{\leaders\hbox to 1em{\hss.\hss}\hfill}
%	use this to write file with table of contents
\def\writetoc{\immediate\openout\tfile=toc.tmp 
   \def\writetoca##1{{\edef\next{\write\tfile{\noindent ##1 
   \string\leaderfill {\noexpand\number\pageno} \par}}\next}}}
%       and this lists table of contents on second pass
%\def\listtoc{\centerline{\bf Contents}\nobreak\medskip{\baselineskip=12pt
%
%  1/98 "toc.tex" has been replaced by "toc.tmp". To get
%  table of content, insert the lines \listtoc and \writetoc, tex twice.
%
% \parskip=0pt\catcode`\@=11 \input toc.tex \catcode`\@=12 \bigbreak\bigskip}}
%

\catcode`\@=12 % at signs are no longer letters
%
%	Unpleasantness in calling in abstract and title fonts
\edef\tfontsize{\ifx\answ\bigans scaled\magstep3\else scaled\magstep4\fi}
\font\titlerm=cmr10 \tfontsize \font\titlerms=cmr7 \tfontsize
\font\titlermss=cmr5 \tfontsize \font\titlei=cmmi10 \tfontsize
\font\titleis=cmmi7 \tfontsize \font\titleiss=cmmi5 \tfontsize
\font\titlesy=cmsy10 \tfontsize \font\titlesys=cmsy7 \tfontsize
\font\titlesyss=cmsy5 \tfontsize \font\titleit=cmti10 \tfontsize
\skewchar\titlei='177 \skewchar\titleis='177 \skewchar\titleiss='177
\skewchar\titlesy='60 \skewchar\titlesys='60 \skewchar\titlesyss='60
\def\titlefont{\def\rm{\fam0\titlerm}% switch to title font
\textfont0=\titlerm \scriptfont0=\titlerms \scriptscriptfont0=\titlermss
\textfont1=\titlei \scriptfont1=\titleis \scriptscriptfont1=\titleiss
\textfont2=\titlesy \scriptfont2=\titlesys \scriptscriptfont2=\titlesyss
\textfont\itfam=\titleit \def\it{\fam\itfam\titleit}\rm}
 \ifx\answ\bigans\else scaled\magstep1\fi
\ifx\answ\bigans\def\abstractfont{\tenpoint}\else
\font\abssl=cmsl10 scaled \magstep1
\font\absrm=cmr10 scaled\magstep1 \font\absrms=cmr7 scaled\magstep1
\font\absrmss=cmr5 scaled\magstep1 \font\absi=cmmi10 scaled\magstep1
\font\absis=cmmi7 scaled\magstep1 \font\absiss=cmmi5 scaled\magstep1
\font\abssy=cmsy10 scaled\magstep1 \font\abssys=cmsy7 scaled\magstep1
\font\abssyss=cmsy5 scaled\magstep1 \font\absbf=cmbx10 scaled\magstep1
\skewchar\absi='177 \skewchar\absis='177 \skewchar\absiss='177
\skewchar\abssy='60 \skewchar\abssys='60 \skewchar\abssyss='60
\def\abstractfont{\def\rm{\fam0\absrm}% switch to abstract font
\textfont0=\absrm \scriptfont0=\absrms \scriptscriptfont0=\absrmss
\textfont1=\absi \scriptfont1=\absis \scriptscriptfont1=\absiss
\textfont2=\abssy \scriptfont2=\abssys \scriptscriptfont2=\abssyss
\textfont\itfam=\bigit \def\it{\fam\itfam\bigit}\def\footnotefont{\tenpoint}%
\textfont\slfam=\abssl \def\sl{\fam\slfam\abssl}%
\textfont\bffam=\absbf \def\bf{\fam\bffam\absbf}\rm}\fi
\def\tenpoint{\def\rm{\fam0\tenrm}% switch back to 10-point type
\textfont0=\tenrm \scriptfont0=\sevenrm \scriptscriptfont0=\fiverm
\textfont1=\teni  \scriptfont1=\seveni  \scriptscriptfont1=\fivei
\textfont2=\tensy \scriptfont2=\sevensy \scriptscriptfont2=\fivesy
\textfont\itfam=\tenit \def\it{\fam\itfam\tenit}\def\footnotefont{\ninepoint}%
\textfont\bffam=\tenbf \def\bf{\fam\bffam\tenbf}\def\sl{\fam\slfam\tensl}\rm}
\font\ninerm=cmr9 \font\sixrm=cmr6 \font\ninei=cmmi9 \font\sixi=cmmi6 
\font\ninesy=cmsy9 \font\sixsy=cmsy6 \font\ninebf=cmbx9 
\font\nineit=cmti9 \font\ninesl=cmsl9 \skewchar\ninei='177
\skewchar\sixi='177 \skewchar\ninesy='60 \skewchar\sixsy='60 
\def\ninepoint{\def\rm{\fam0\ninerm}% switch to footnote font
\textfont0=\ninerm \scriptfont0=\sixrm \scriptscriptfont0=\fiverm
\textfont1=\ninei \scriptfont1=\sixi \scriptscriptfont1=\fivei
\textfont2=\ninesy \scriptfont2=\sixsy \scriptscriptfont2=\fivesy
\textfont\itfam=\ninei \def\it{\fam\itfam\nineit}\def\sl{\fam\slfam\ninesl}%
\textfont\bffam=\ninebf \def\bf{\fam\bffam\ninebf}\rm} 
%
%---------------------------------------------------------------------
%

\hyphenation{anom-aly anom-alies coun-ter-term coun-ter-terms}
\def\inv{^{\raise.15ex\hbox{${\scriptscriptstyle -}$}\kern-.05em 1}}

\def\Dsl{\,\raise.15ex\hbox{/}\mkern-13.5mu D} %this one can be subscripted
\def\dsl{\raise.15ex\hbox{/}\kern-.57em\partial}

\font\bigit=cmti10 scaled \magstep1
 %pound sterling
\def\lspace{\ifx\answ\bigans{}\else\qquad\fi}
\def\lbspace{\ifx\answ\bigans{}\else\hskip-.2in\fi} % $$\lbspace...$$
\def\boxeqn#1{\vcenter{\vbox{\hrule\hbox{\vrule\kern3pt\vbox{\kern3pt
	\hbox{${\displaystyle #1}$}\kern3pt}\kern3pt\vrule}\hrule}}}
\def\mbox#1#2{\vcenter{\hrule \hbox{\vrule height#2in
		\kern#1in \vrule} \hrule}}  %e.g. \mbox{.1}{.1}
%	matters of taste
%\def\tilde{\widetilde} \def\bar{\overline} \def\hat{\widehat}
%
% some sample definitions
  %     curly letters
  \def\CF{{\cal F}}

\def\darr#1{\raise1.5ex\hbox{$\leftrightarrow$}\mkern-16.5mu #1}
 %pound sterling

 %puts a small half in a displayed eqn
\def\roughly#1{\raise.3ex\hbox{$#1$\kern-.75em\lower1ex\hbox{$\sim$}}}

\def\frac#1#2{{#1\over#2}}

\def\journal#1&#2(#3){\unskip, #1~\bf #2 \rm(19#3) }
\def\andjournal#1&#2(#3){\sl #1~\bf #2 \rm (19#3) }

\def\bra#1{\left\langle #1\right|}
\def\ket#1{\left| #1\right\rangle}

\def\One{{1\hskip -3pt {\rm l}}}
\catcode`\@=11\def\slash#1{\mathord{\mathpalette\c@ncel{#1}}}
\overfullrule=0pt
\def\steepslash{\c@ncel}
\def\frac#1#2{{#1\over #2}}

\def\:{\!:\!}
\def\inbar{\,\vrule height1.5ex width.4pt depth0pt}
\def\IQ{\relax\,\hbox{$\inbar\kern-.3em{\rm Q}$}}
\def\IB{\relax{\rm I\kern-.18em B}}
\def\IC{\relax\hbox{$\inbar\kern-.3em{\rm C}$}}
\def\IP{\relax{\rm I\kern-.18em P}}
\def\IR{\relax{\rm I\kern-.18em R}}
\def\ZZ{\relax\ifmmode\mathchoice
{\hbox{Z\kern-.4em Z}}{\hbox{Z\kern-.4em Z}}
{\lower.9pt\hbox{Z\kern-.4em Z}}
{\lower1.2pt\hbox{Z\kern-.4em Z}}\else{Z\kern-.4em Z}\fi}

\catcode`\@=12

%                      Zeitschriften:
\def\npb#1(#2)#3{{ Nucl. Phys. }{B#1} (#2) #3}
\def\plb#1(#2)#3{{ Phys. Lett. }{#1B} (#2) #3}
\def\pla#1(#2)#3{{ Phys. Lett. }{#1A} (#2) #3}
\def\prl#1(#2)#3{{ Phys. Rev. Lett. }{#1} (#2) #3}
\def\mpla#1(#2)#3{{ Mod. Phys. Lett. }{A#1} (#2) #3}
\def\ijmpa#1(#2)#3{{ Int. J. Mod. Phys. }{A#1} (#2) #3}
\def\cmp#1(#2)#3{{ Comm. Math. Phys. }{#1} (#2) #3}
\def\cqg#1(#2)#3{{ Class. Quantum Grav. }{#1} (#2) #3}
\def\jmp#1(#2)#3{{ J. Math. Phys. }{#1} (#2) #3}
\def\anp#1(#2)#3{{ Ann. Phys. }{#1} (#2) #3}
\def\prd#1(#2)#3{{ Phys. Rev. } {D{#1}} (#2) #3}
\def\ptp#1(#2)#3{{ Progr. Theor. Phys. }{#1} (#2) #3}
\def\aom#1(#2)#3{{ Ann. Math. }{#1} (#2) #3}

\def\br{\buildrel}
\def\bra{\langle}
\def\ket{\rangle}

\def\C{{\bf C}}

\def\O{{\bf O}}
\def\P{{\bf P}}

\def\Z{{\bf Z}}

\def\cO{{\cal O}}

\def\cQ{{\cal Q}}
\def\cR{{\cal R}}

\def\cT{{\cal T}}
\def\cU{{\cal U}}

\def\cW{{\cal W}}

\def\cY{{\cal Y}}

\def\cicy#1(#2|#3)#4{\left(\matrix{#2}\right|\!\!
                     \left|\matrix{#3}\right)^{{#4}}_{#1}}

\def\lra{\longrightarrow}
\def\lla{\longleftarrow}
\def\ra{\rightarrow}

\def\bs{\bigskip}

\def\Box{{\,\lower0.9pt\vbox{\hrule 
\hbox{\vrule height 0.2 cm \hskip 0.2 cm  
\vrule height 0.2 cm}\hrule}\,}}

\global\newcount\thmno \global\thmno=0
\def\definition#1{\global\advance\thmno by1
\bigskip\noindent{\bf Definition \secsym\the\thmno. }{\it #1}
\par\nobreak\medskip\nobreak}
\def\question#1{\global\advance\thmno by1
\bigskip\noindent{\bf Question \secsym\the\thmno. }{\it #1}
\par\nobreak\medskip\nobreak}
\def\theorem#1{\global\advance\thmno by1
\bigskip\noindent{\bf Theorem \secsym\the\thmno. }{\it #1}
\par\nobreak\medskip\nobreak}
\def\proposition#1{\global\advance\thmno by1
\bigskip\noindent{\bf Proposition \secsym\the\thmno. }{\it #1}
\par\nobreak\medskip\nobreak}
\def\corollary#1{\global\advance\thmno by1
\bigskip\noindent{\bf Corollary \secsym\the\thmno. }{\it #1}
\par\nobreak\medskip\nobreak}
\def\lemma#1{\global\advance\thmno by1
\bigskip\noindent{\bf Lemma \secsym\the\thmno. }{\it #1}
\par\nobreak\medskip\nobreak}
\def\conjecture#1{\global\advance\thmno by1
\bigskip\noindent{\bf Conjecture \secsym\the\thmno. }{\it #1}
\par\nobreak\medskip\nobreak}
\def\exercise#1{\global\advance\thmno by1
\bigskip\noindent{\bf Exercise \secsym\the\thmno. }{\it #1}
\par\nobreak\medskip\nobreak}
\def\remark#1{\global\advance\thmno by1
\bigskip\noindent{\bf Remark \secsym\the\thmno. }{\it #1}
\par\nobreak\medskip\nobreak}
\def\problem#1{\global\advance\thmno by1
\bigskip\noindent{\bf Problem \secsym\the\thmno. }{\it #1}
\par\nobreak\medskip\nobreak}
\def\others#1#2{\global\advance\thmno by1
\bigskip\noindent{\bf #1 \secsym\the\thmno. }{\it #2}
\par\nobreak\medskip\nobreak}
\def\proof{\noindent Proof: }

\def\thmlab#1{\xdef #1{\secsym\the\thmno}\writedef{#1\leftbracket#1}\wrlabeL{#1=#1}}
%
% redefine \newsec so that all \thmno set to zero in a new section
%
\def\newsec#1{\global\advance\secno by1\message{(\the\secno. #1)}
%\ifx\answ\bigans \vfill\eject \else \bigbreak\bigskip \fi  %if desired
\global\subsecno=0\thmno=0\eqnres@t\noindent{\bf\the\secno. #1}
\writetoca{{\secsym} {#1}}\par\nobreak\medskip\nobreak}
\def\eqnres@t{\xdef\secsym{\the\secno.}\global\meqno=1\bigbreak\bigskip}
\def\sequentialequations{\def\eqnres@t{\bigbreak}}\xdef\secsym{}
%

%

%%%%%%%%%%%%%%%%%%%%%%
%  REFERENCES        %
%%%%%%%%%%%%%%%%%%%%%%

\ref\AP{C. Allday, V. Puppe,{\it Cohomology methods in transformation groups},
Cambridge University Press, 1993.}
\ref\AB{M. Atiyah and R. Bott, {\it The moment map and
equivariant cohomology}, Topology 23 (1984) 1-28.} 
\ref\Audin{M. Audin, {\it The topology of torus actions
on symplectic manifolds}, Prog. in Math. 93,
Birkh\"auser 1991.}
\ref\Batyrev{V. Batyrev, {\it Dual polyhedra and mirror symmetry
for Calabi-Yau hypersurfaces in toric varieties},
J. Alg. Geom. 3 (1994) 493-535.}
\ref\BB{V. Batyrev and L. Borisov, {\it 
On Calabi-Yau complete intersections in toric varieties},
alg-geom/9412017.}
\ref\BehrendFentachi{K. Behrend and B. Fantechi, 
{\it The intrinsic normal cone}, Invent. Math. 128 (1997) 45-88.}
\ref\BV{N. Berline and M. Vergne, {\it Classes caracte\'eristiques
equivariantes, Formula de localisation en cohomologie
equivariante.} 1982.}
\ref\Be{A. Bertram, {\it Towards a Schubert calculus for maps from 
Riemann surface to Grassmannian}, Inter. J. Math., V.5, No.6 (1994) 811-825.}
\ref\CDGP{P. Candelas, X. de la Ossa, P. Green, and
L. Parkes, {\it A pair of Calabi-Yau manifolds as an
exactly soluble superconformal theory},
Nucl. Phys. B359 (1991) 21-74.}
\ref\CF{I. Ciocan-Fontanine, {\it Quantum cohomology of flag manifolds},
IMRN, 1995, No. 6, 263-277.}
\ref\Cox{D. Cox, {\it
Homogeneous coordinate ring of a toric variety}, J. Algebraic Geom. 4 (1995) 17-50.} 
\ref\EG{D. Edidin and W. Graham,{\it 
 Equivariant intersection theory}, alg-geom/9609018.}
\ref\Fulton{W. Fulton, {\it Intersection Theory}, Springer-Verlag,
2nd Ed., 1998.}
\ref\FO{K. Fukaya and K. Ono, {\it
Arnold conjecture and Gromov-Witten invariants,}
preprint, 1996.}
\ref\GKZ{I. Gel'fand, M. Kapranov and A. Zelevinsky,
{\it Hypergeometric functions and toral manifolds},
Funct. Anal. Appl. 23 (1989) 94-106.}
\ref\GiventalI{A. Givental, {\it A mirror theorem for
toric complete intersections}, alg-geom/9701016.}
\ref\GP{T. Graber and R. Pandharipande, {\it Localization of
virtual classes}, alg-geom/9708001.}
\ref\GKM{M. Goresky, R. Kottwitz and R. MacPherson,
{\it Equivariant cohomology, Koszul duality and the
localization theorem}, Invent. Math. 131 (1998) 25-83.}
\ref\GZ{V. Guillemin and C. Zara, {\it Equivariant de Rham
Theory and Graphs}, math.DG/9808135.}
\ref\Hirzebruch{F. Hirzebruch, {\it Topological methods in
algebraic geometry}, Springer-Verlag, Berlin 1995, 3rd Ed.}
\ref\HKTY{S. Hosono, A. Klemm,
           S. Theisen and S.T. Yau, {\it Mirror symmetry, mirror map
           and applications to complete intersection Calabi-Yau
           spaces}, hep-th/9406055, Nucl. Phys. B433 (1995) 501-554.}
\ref\HLYII{S. Hosono, B.H. Lian and S.T. Yau, {\it GKZ-generalized
hypergeometric systems and mirror symmetry of Calabi-Yau hypersurfaces},
alg-geom/9511001, Commun. Math. Phys. 182 (1996) 535-578.}
\ref\HLY{S. Hosono, B.H. Lian and S.T. Yau, {\it Maximal
Degeneracy Points of GKZ Systems}, alg-geom/9603014,
Journ. Am. Math. Soc., Vol. 10, No. 2 (1997) 427-443.}
\ref\Jaczewski{K. Jaczewski, {\it Generalized Euler sequence
and toric varieties}, Contemporary Math. Vol 162 (1994) 227-247.}
\ref\Kr{A. Kresch,{\it Cycle groups for Artin Stacks,} math. AG/9810166.}
\ref\Kontsevich{M. Kontsevich,
{\it Enumeration of rational curves via torus actions,}
in {\it The Moduli Space of Curves}, ed. by
R. Dijkgraaf, C. Faber, G. van der Geer, Progress in Math.
vol. 129, Birkh\"auser, 1995, 335--368.}
\ref\LiTianII{ J. Li and G. Tian, {\it
Virtual moduli cycle and 
Gromov-Witten invariants of algebraic varieties},
J. of Amer. math. Soc. 11, no. 1, (1998) 119-174.}
\ref\LiTianV{J. Li and G. Tian, {\it A brief tour of GW invariants},
to appear in Surveys in Diff. Geom. Vol. 6 (1999).}
\ref\LLYI{B. Lian, K. Liu and S.T. Yau, {\it Mirror Principle I},
Asian J. Math. Vol. 1, No. 4 (1997) 729-763.}
\ref\LLYII{B. Lian, K. Liu and S.T. Yau, {\it Mirror Principle II},
To appear in Asian J. Math. (1999).}
\ref\MP{D. Morrison and R. Plesser, {\it
Summing the instantons: quantum cohomology and mirror symmetry
in toric varieties}, alg-geom/9412236.}
\ref\Musson{I. Musson, {\it Differential operators on toric
varieties}, JPAA 95 (1994) 303-315.}
\ref\Oda{T. Oda, {\it Convex Bodies and Algebraic Geometry},
Series of Modern Surveys in Mathematics 15,
Springer Verlag (1985).}
\ref\Ruan{Y.B. Ruan, {\it Virtual neighborhoods 
and pseudo-holomorphic curves}, math.AG/9611021.}
\ref\RuanTian{Y.B. Ruan and G. Tian, {\it A mathematical
theory of quantum cohomology}, Journ. Diff. Geom. Vol. 42,
No. 2 (1995) 259-367.}
\ref\SiII{B. Siebert, {\it
Gromov-Witten invariants for general symplectic manifolds,}
alg-geom/9608005.}
\ref\So{F. Sottile, {\it Real rational curves in Grassmannians }, Math.AG/9904167.}  
\ref\SS{F. Sottile, B. Sturmfels, {\it A sagbi basis for the quantum 
Grassmannian }, Math.AG/9908016.} 
\ref\Stromme{S.A. Stromme, {\it On parametrized rational curves in Grassmann varieties},
in {\it Space curves}, Ed. Ghione et al. LNM 1266, Springer-Verlag.}
\ref\VistoliI{A. Vistoli, {\it Intersection theory on algebraic stacks
and their moduli}, Inv. Math. 97 (1989) 613-670.}
\ref\VistoliII{A. Vistoli, {\it Equivariant Grothendieck groups and
equivariant Chow groups in Classification of irregular varieties},
LNM 1515 (1992) 112-133.}
\ref\Witten{E. Witten, {\it Phases of N=2 theories in two
dimension}, hep-th/9301042.}

\Title{}{Mirror Principle III} 

\centerline{
Bong H. Lian,$^{1}$\footnote{}{$^1$~~Department of Mathematics,
Brandeis University, Waltham, MA 02154.}
Kefeng Liu,$^{2}$\footnote{}{$^2$~~Department of Mathematics,
Stanford University, Stanford, CA 94305.}
and Shing-Tung Yau$^3$\footnote{}{$^3$~~Department of Mathematics,
Harvard University, Cambridge, MA 02138.} }

\vskip .2in

Abstract. 
We generalize the theorems in {\it Mirror Principle I} and {\it II}
to the case of general projective manifolds without the
convexity assumption. We also apply the results to
balloon manifolds, and generalize to higher genus.

%\listtoc %input toc.tmp
\writetoc %create toc.tmp

\Date{}

\newsec{Introduction}

The present paper is a sequel to {\it Mirror Principle I} and {\it II} \LLYI
\LLYII. For motivations and the main ideas of mirror principle,
we refer the reader to the introductions of these two papers.
 
Let $X$ be a projective manifold, and $d\in A_1(X)$.
Let $M_{0,k}(d,X)$ denote the moduli space of
$k$-pointed,  genus 0, degree $d$, stable maps
$(C,f,x_1,..,x_k)$ with target $X$
\Kontsevich. 
{\it Note that our notation is without the bar.}
By the construction of \LiTianII
(see also \BehrendFentachi\FO),
each nonempty  
$M_{0,k}(d,X)$ admits a homology class $LT_{0,k}(d,X)$ 
of dimension $dim~X+\bra c_1(X),d\ket +k-3$. This
class plays the role of the fundamental class in
topology, hence
$LT_{0,k}(d,X)$ is
called the virtual fundamental class.
For background on this, we recommend \LiTianV.

Let $V$ be a convex vector bundle on $X$.
(ie. $H^1(\P^1,f^*V)=0$ for every holomorphic
map $f:\P^1\ra X$.) Then $V$ induces on each
$M_{0,k}(d,X)$ a vector bundle $V_d$,
with fiber at
$(C,f,x_1,..,x_k)$ given by the section space $H^0(C,f^*V)$.
Let $b$ be any multiplicative characteristic class
\Hirzebruch.
(ie. if $0\ra E'\ra E\ra E''\ra 0$ is an exact sequence
of vector bundles,
then $b(E)=b(E')b(E'')$.) 
The problem we study here is to understand
the intersection numbers
$$K_d:=\int_{LT_{0,0}(d,X)} b(V_d)$$
and their generating function:
$$\Phi(t):=\sum K_d~ e^{d\cdot t}.$$
There is a similar and equally important problem if one
starts from a concave vector bundle $V$ \LLYI.
(ie. $H^0(\P^1,f^*V)=0$ for every holomorphic
map $f:\P^1\ra X$.) More generally, $V$ can
be a direct sum of a convex and a concave bundle.
Important progress made on these problems has come from
mirror symmetry. All of it
seems to point toward the following
general phenomenon \CDGP, which we call
{\it the Mirror Principle}. Roughly,
it says that there are functional identities
which can be used to either constrain or to
compute the $K_d$
often in terms of certain explicit
special functions, loosely called generalized
hypergeometric functions. In this paper, we
generalize this principle to include all
projective manifolds. We apply this theory
to compute the multiplicative classes $b(V_d)$ for vector bundles
on balloon manifolds. The answer is in terms of
certain universal virtual classes which are independent of $V, b$.

When $X$ is a toric manifold,
$b$ is the Euler class, and $V$ is a sum
of line bundles, there is a general formula
derived in \HKTY\HLY~ based on mirror symmetry, giving
$\Phi(t)$ in terms of generalized hypergeometric functions \GKZ.
Similar functions were studied \GiventalI~
in equivariant quantum cohomology theory
based on a series of axioms. For further background, 
see introduction of \LLYI.

{\bf Acknowledgements.} We thank Y. Hu, C.H. Liu, and G. Tian,
for numerous helpful discussions with us during the course
of this project. We owe special thanks to J. Li for
patiently explaining to us his joint work with Tian,
for tirelessly providing a lot of technical assistance, 
and for proofreading a substantial part of this manuscript.
B.H.L.'s research
is supported by NSF grant DMS-9619884.
K.L.'s research is supported by NSF grant
DMS-9803234 and the Terman fellowship.
S.T.Y.'s research is supported by DOE grant
DE-FG02-88ER25065 and NSF grant DMS-9803347.

\subsec{Outline}

In section 2, we do the necessary preparation to
set up the version of  localization theorem we need.
This is a (functorial localization) formula which translates
a commutative square diagram into
a relation between localizations on two 
$T$ spaces related by an equivariant map.

We do basically three things in section 3. 
After we introduced the necessary notations,
first we apply functorial localization
to stable map moduli spaces.
Second, 
we prove one of the main results of this paper:
Theorem 3.6, which translate structure of fixed points
on stable map moduli into an algebraic identity on the
homology of a projective manifold (with or without $T$ action).
This motivates the  notion of Euler data and
Euler series. These are essentially solutions to
the algebraic identity just mentioned. 
Third, we prove the main Theorems 3.12-3.13 which relate the generating
functions $\Phi(t)$ with an Euler series $A(t)$ arising from
induced bundles on stable map moduli. 

In section 4, we specialize results in section 3 to
balloon manifolds, and introduce the notion of
linking. The main theorems here are 4.5 and 4.7.
The first of these gives a description of
an essential polar term of $A(t)$ upon localizing
at a fixed point in $X$. The second theorem
gives a sufficient condition for computing $A(t)$
in terms of certain universal virtual classes
on stable map moduli. We then specialize this to the
case when $b_T$ is the Euler class or the Chern polynomial.

In section 5, we explain some other ways to compute $A(t)$,
first by relaxing those sufficient conditions, then
by finding an explicit closed formula for
those universal virtual classes above by
using an equivariant short exact sequence for the tangent bundle.
This includes toric manifolds
as a special case. We then formulate
an inductive method for computing $A(t)$ in
full generality for any balloon manifold.
Next, we discuss a method in which 
functorial localization is used to study $A(t)$
via a resolution of the image of the collapsing map.
In certain cases, this resolution can be
described quite explicitly. Finally, we discuss
a generalization of mirror principle to higher genus.

\newsec{Set-up}

Basic references: on intersection theory 
on algebraic schemes and stacks, we use \Fulton\VistoliI;
on the virtual classes,
we follow \LiTianII;
on their equivariant counterparts,
see \AP\AB\BV\Kr\EG\GP\VistoliII.

$T$ denotes an algebraic torus.
$T$-equivariant Chow groups (homology) with complex coefficients
are denoted by $A^T_*(\cdot)$.
$T$-equivariant operational Chow groups (cohomology) with complex coefficients
are denoted by $A_T^*(\cdot)$. 
For $c\in A^p_T(X)$, and $\beta\in A_q^T(X)$, we denote
by $c\cap \beta=\beta\cap c$ the image of $c\otimes\beta$
under the canonical homomorphisms
$$A^p_T(X)\otimes A_q^T(X)\ra A_{q-p}^T(X).$$
The product on $A_T^*(X)$ is denoted by $a\cdot b$.
The homomorphisms $\cap$ define an $A^*_T(X)$-module structure
on the homology $A_*^T(X)$.
When $X$ is nonsingular, there is a
compatible intersection product on $A^T_*(X)$
which we denote by  $\beta\cdot \gamma$.

Given a $T$-equivariant (proper or flat) map $f:X\ra Y$, we denote by
$$f_*:A^T_*(X)\ra A^T_*(Y),~~~~f^*:A^T_*(Y)\ra A^T_*(X)$$
the equivariant (proper) pushforward and (flat) pullback;
the notations $f^*$ and $f_*$ are also used for pullback
and (flat) pushforward on cohomology.
All maps used here will be assumed proper.
A formula often used is the projection formula:
$$f_* (f^* c\cap\beta)=c\cap f_*(\beta)$$
for cohomology class $c$ on $Y$ and homology class $\beta$ on $X$.
Note that both $A^T_*(X)$ and $A_T^*(X)$ are modules over
the algebra $A^*_T(pt)=\C[\cT^*]$,
where $\cT^*$ is the dual of the Lie algebra of $T$,
and the homomorphisms $f_*,f^*$ are module homomorphisms.
We often extend these homomorphisms over the field $\C(\cT^*)$
without explicitly saying so.
%The $T$ actions on schemes and stacks considered here
%are assumed effective.
%If both $X,Y$ are nonsingular, then 
%pushforward also exists on cohomology such that
%$$f_*:A_T(X)\ra A_T(Y),~~~~f_*(c)\cap\beta=f_*(c\cap f^*\beta),~
%\beta\in A^T(Y).$$
Finally, suppose we have a fiber square
$${\matrix{
F & {\br i'\over\longrightarrow}& M\cr
p\downarrow &  &\downarrow q\cr
X & {\br i\over\longrightarrow} & Y
}}
$$
where $i$ is a regular embedding of codimension $d$, then we have
$$p_*i^!\beta=i^* q_*\beta$$
for any homology class $\beta$ on $M$.
Here $i^!:A^T_*(M)\ra A^T_{*-d}(F)$
is the refined Gysin homomorphism.

\subsec{Functorial localization}

Let $X$ be an algebraic scheme with a $T$ action and
equipped with a suitable perfect
obstruction theory (see \LiTianII\GP). 
Let $F_r$ denote the fixed point components in $X$.
Let $[X]^{vir},[F_r]^{vir}$ be the equivariant virtual classes 
of $X$ and the $F_r$. Then by \GP, 
$$
[X]^{vir}=\sum_r {i_r}_* {[F_r]^{vir}\over e_T(F_r/X)}
$$
where $i_r:F_r\ra X$ are the inclusions, and $e_T(F_r/X)$
the equivariant Euler class of the virtual normal bundle of $F_r\subset X$.
Then for any cohomology class $c$ on $X$, we have
\eqn\dumb{
c\cap [X]^{vir}=
\sum_r {i_r}_* {i^*_r c\cap[F_r]^{vir}\over e_T(F_r/X)}.
}

Throughout this subsection, let
$$f:\ \ X \rightarrow Y$$
be an equivariant map with $Y$ smooth. 
Let $E$ be a fixed point component in $Y$, and
let $F$ be the fixed points in $f^{-1}(E)$. Let $g$ be the restriction
of $f$ to $F$, and $j_E:E\ra Y$, $i_F:F\ra X$ be the inclusion maps.
Thus we have the commutative diagram:
\eqn\CommSquare{\matrix{
F & {\br i_F\over \longrightarrow} & X\cr
g\downarrow &  & \downarrow f\cr
E & {\br j_E\over \longrightarrow} & Y.}
}
Then we have the following {\it functorial localization formula}.

\lemma{ Given a cohomology class $\omega \in A^*_T(X)$, we have the
equality on $E$:
$$
{j_E^*f_*(\omega\cap [X]^{vir})\over e_T(E/Y)}
=g_*\left({i_F^*\omega\cap[F]^{vir}\over e_T(F/X)}\right).
$$
}
\proof Applying \dumb~ to
the class $c=\omega\cdot f^*{j_E}_*1$ on $X$, we get
$$\omega\cdot f^*{j_E}_*1\cap [X]^{vir}={i_F}_*\left(
{i_F^*(\omega
\cdot f^*{j_E}_*1)\cap [F]^{vir}\over e_T(F/X)}\right).$$ 
Note that the contributions from fixed
components other than $F$ vanish. Applying
$f_*$ to both
sides, we get
$$f_*(\omega\cap [X]^{vir})\cap {j_E}_*1=f_*{i_F}_*
\left({i_F^*(\omega\cdot
f^*{j_E}_*1)\cap [F]^{vir}\over e_T(F/X)}\right).$$
Now $f\circ i_F=j_E\circ g$ which, implies
$$f_*{i_F}_*={j_E}_*g_*,\ \ i_F^* f^*=g^*j_E^* .$$ 
Thus we get
$$f_*(\omega\cap [X]^{vir})\cap {j_E}_*1={j_E}_*g_*
\left({i_F^*\omega\cdot g^*e_T(E/Y)\cap [F]^{vir}\over e_T(F/X)}\right).$$
Applying $j_E^*$ to both sides here, we get
$$\eqalign{
j_E^*f_*(\omega\cap[X]^{vir})\cap e_T(E/Y)&=e_T(E/Y)\cap
g_*\left({i_F^*\omega\cdot g^*e_T(E/Y)\cap[F]^{vir}\over e_T(F/X)}\right)\cr
&=e_T(E/Y)^2\cap g_*\left({i_F^*\omega\cap[F]^{vir}\over e_T(F/X)}\right).
}$$
Since $e_T(E/Y)$ is invertible, our assertion follows.
$\Box$

Note that  if $F$ has more than one component, then
the right hand side of the formula above becomes a sum
over those components in an obvious way.

\corollary{Let $Y'$ be a $T$-invariant submanifold
of $Y$, $f':X'=f^{-1}(Y')\ra Y'$ be the restriction of $f:X\ra Y$
to the substack $X'$,
and $j:Y'\ra Y$, $i:X'\ra X$ be the inclusions. Then
for any $\omega \in A^*_T(X)$, we have 
$${j^*f_*(\omega\cap [X]^{vir})\over e_T(Y'/Y)}
=f'_*\left({i_F^*\omega\cap[X']^{vir}\over e_T(X'/X)}\right).
$$
}
\proof
Let $E$ be any fixed point component of $Y$ contained in $Y'$,
and $F$ be the fixed points in $f^{-1}(E)$, as in the
preceding lemma. Then we have the commutative diagram
\eqn\dumb{\matrix{
F & {\br i_F'\over \longrightarrow} & X'
 & {\br i\over \longrightarrow} & X\cr
g\downarrow &   &  f'\downarrow &   & f\downarrow\cr
E & {\br j_E'\over \longrightarrow} & Y'
 & {\br j\over \longrightarrow} & Y.}
}
We will show that
$$(*)~~~~{j_E'}^*\left({j^*f_*(\omega\cap [X]^{vir})\over e_T(Y'/Y)}\right)
={j_E'}^*f'_*\left({i_F^*\omega\cap[X']^{vir}\over e_T(X'/X)}\right).
$$
Then our assertion follows from the localization theorem.

Put $j_E:=j\circ j_E'$, $i_F:=i\circ i_F'$. The left hand side of (*) is
$$\eqalign{
{ {j_E'}^*j^*f_*(\omega\cap[X]^{vir})\over
{j_E'}^* e_T(Y'/Y)}
&=e_T(E/Y')\cap {j_E^*f_*(\omega\cap[X]^{vir})\over e_T(E/Y)}\cr
&=e_T(E/Y')\cap 
g_*\left({i_F^*\omega\cap[F]^{vir}\over e_T(F/X)}\right)
~~~~~(preceding~ lemma).
}$$
Now apply the left hand square in \dumb~ and the preceding lemma again
to the class ${i^*\omega\over e_T(X'/X)}$ on $X'$.
Then the right hand side of (*) becomes
$$\eqalign{
{j_E'}^*f'_*\left({i_F^*\omega\over e_T(X'/X)}\cap[X']^{vir}\right)
&=e_T(E/Y')\cap
g_*\left({ {i_F'}^*{i^*\omega\over e_T(X'/X)}\cap[F]^{vir}\over e_T(F/X')}\right)\cr
&=e_T(E/Y')\cap 
g_*\left({i_F^*\omega\cap[F]^{vir}\over e_T(F/X)}\right).
}
$$
This proves (*). $\Box$

\newsec{General Projective $T$-manifolds}

Let $X$ be a projective $T$-manifold.
Let $M_d(X)$ be the degree $(1,d)$, arithmetic genus zero, 0-pointed, stable map
moduli stack with target $\P^1\times X$. 
The standard $\C^\times$ action 
on $\P^1$ together with the $T$ action on $X$ induces
a $G=\C^\times\times T$ action on $M_d(X)$. 
Let $LT_d(X)\in A_*^G(M_d(X))$
be the virtual class of this moduli stack. This is
an equivariant homology class of dimension $\bra c_1(X),d\ket+dim~X$.

The $\C^\times$
fixed point components $F_r$,
labelled by $0\preceq r\preceq d$, in $M_d(X)$
can be described as follows (see \LLYII).
Let $F_r$ be the substack 
$$F_r:= M_{0,1}(r, X)\times_X M_{0,1}(d-r, X)$$ 
obtained from gluing the two one pointed moduli stacks. 
More precisely, consider the map
$$e^X_r\times  e^X_{d-r}: 
M_{0,1}(r, X)\times M_{0,1}(d-r, X)\rightarrow X\times X$$ 
given by evaluations at the corresponding marked points; and
$$\Delta:\ X\rightarrow X\times X$$ 
the diagonal map. Then we have
$$F_r=(e^X_r\times e^X_{d-r})^{-1}\Delta(X).$$
Note that $F_d=M_{0,1}(d,X)=F_0$ by convention,
but $F_0$ and $F_d$ will be embedded into $M_d(X)$
in two different ways.
The  $F_r$ can be identified with a $\C^\times$ fixed point
component of $M_d(X)$ as follows. 
Consider the case $r\neq0,d$ first. Given a pair
$(C_1,f_1,x_1)\times (C_2,f_2,x_2)$ in $F_r$, we
get a new curve $C$ by gluing
$C_1,C_2$ to $\P^1$ with
$x_1,x_2$ glued to $0,\infty\in\P^1$ respectively.
The new curve $C$ is mapped into $\P^1\times X$ as follows.
Map $\P^1\subset C$ identically onto $\P^1$, and
collapse $C_1,C_2$ to $0,\infty$ respectively; 
then map $C_1,C_2$ into $X$ with $f_1,f_2$ respectively,
and collapse the $\P^1$ to $f(x_1)=f(x_2)$. 
This defines a stable map
$(C,f)$ in $M_d(X)$. 
For $r=d$, we glue $(C_1,f_1,x_1)$ to
$\P^1$ at $x_1$ and $0$. For $r=0$, we glue
$(C_2,f_2,x_2)$ to $\P^1$ at $x_2$ and $\infty$.

{\it Notations:}
\item{(i)} We identify $F_r$
as a substack of $M_d(X)$ as above, and let
$$i_r:F_r\ra M_d(X)$$
denote  the inclusion map. 

\item{(ii)}
We have evaluation maps
$$e^X:F_r\ra X,$$
which sends a pair in $F_r$ to the value
at the common marked point. While the
notation $e^X$ doesn't reflect the dependence
on $r$, the domain $F_r$ that $e^X$ operates on
will be clear.

\item{(iii)} We have the obvious inclusion
$$\Delta':F_r\subset M_{0,1}(r,X)\times M_{0,1}(d-r,X),$$
and projections 
$$p_0:F_r\ra M_{0,1}(r,X),~~~~p_\infty:F_r\ra M_{0,1}(d-r,X).$$

\item{(iv)}
Let $L_r$ denote the universal line bundle on $M_{0,1}(r,X)$.

\item{(v)} We have the natural forgetting, evaluation,
and projection maps:
$$\eqalign{
&\rho:M_{0,1}(d,X)\ra M_{0,0}(d,X)\cr
&e^X_d:M_{0,1}(d,X)\ra X\cr
&\pi:M_d(X)\ra M_{0,0}(d,X).
}$$
We also have the obvious commutative diagrams
$$\matrix{
M_d(X) &  &   \cr
\pi\downarrow & \nwarrow i_0 &  \cr
M_{0,0}(d,X) & {\br \rho\over\longleftarrow} & M_{0,1}(d,X)
}$$
\eqn\FiberSquare{\matrix{
F_r & {\br \Delta'\over\longrightarrow}& 
M_{0,1}(r,X)\times M_{0,1}(d-r,X)\cr
e^X\downarrow &  &~~~~~~~\downarrow e^X_r\times e^X_{d-r}\cr
X & {\br \Delta\over\longrightarrow} & X\times X
}}
where $\Delta$ is the diagonal map.
Note that we have
a diagram  similar to \FiberSquare~
but with $X$ replaced by $Y$ in the bottom row.
From the fiber square \FiberSquare, we have a refined Gysin homomorphism
$$\Delta^!:A^T_*(
M_{0,1}(r,X)\times M_{0,1}(d-r,X))\ra 
A^T_{*-dim~X}(F_r).
$$
We refer the reader to section 6 \LiTianII~ for the following

\lemma{\LiTianII~
For $r\neq 0,d$,
$[F_r]^{vir}=\Delta^!(LT_{0,1}(r,X)\times LT_{0,1}(d-r,X))$.
}

\item{(vi)} Let $\alpha$ be the weight of
the standard $\C^\times$ action on $\P^1$.
We denote by $A_*^T(X)(\alpha)$ the algebra obtained from
the polynomial algebra
$A_*^T(X)[\alpha]$ by localizing 
with respect to all invertible elements.
If $\beta$ is an element in $A_*^T(X)(\alpha)$,
we let $\overline{\beta}$ be the class obtained from $\beta$
by replacing $\alpha\ra-\alpha$. We also introduce
formal variables $\zeta=(\zeta_1,...,\zeta_m)$ such that
$\bar\zeta_a=-\zeta_a$.
Denote $\cR=\C(\cT^*)[\alpha]$. When a multiplicative class $b_T$,
such as the Chern polynomial $c_T=x^r+x c_1+\cdots+c_r$,
is considered, we must replace the ground field $\C$ by $\C(x)$,
so that $c_T$ takes value in Chow groups with appropriate coefficients.
This change of ground field will be implicit whenever necessary.

\item{(vii)} For each $d$,
let $\varphi: M_d(X)\ra W_d$ be a
$G$-equivariant map into smooth manifold (or orbifold) $W_d$
with the property
{\it that the $\C^\times$ fixed point components in $W_d$ are
$G$-invariant submanifolds $Y_r$ such that
$\varphi^{-1}(Y_r)=F_r$}.

\item{}
The spaces $W_d$ exist but are not unique. Two specific kinds will
be used here. First, choose an 
equivariant projective embedding
$$\tau: X\ra Y=\P^{n_1}\times\cdots\times\P^{n_m}$$
which induces an isomorphism $A^1(X)\cong A^1(Y)$.
Then we have a $G$-equivariant embedding
$$
M_d(X)\ra M_d(Y).
$$
There is  a $G$-equivariant map (see \LLYI~ and references there)
$$
M_d(Y)\ra W_d:=N_{d_1}\times\cdots\times N_{d_m}
$$
where the $N_{d_a}:=\P H^0(\P^1,\cO(d_a))^{n_a+1}
\cong\P^{(n_a+1)d_a+n_a}$, which are
the linear sigma model
for the $\P^{n_i}$.
Thus composing the two maps above, we get
a $G$-equivariant map $\varphi:M_d(X)\ra W_d$.
It is also easy to check that the $\C^\times$
fixed point components in $W_d$ have the desired property.
Second, if
$X$ is a toric variety, then there exist toric varieties $W_d$
\MP~ where $Y_r$ are submanifolds of $X$.
Note that $X$ is contained in the $Y_r=Y$ as a submanifold in the first
kind, while $X$ contains the $Y_r$ as submanifolds in the
second kind.
We will used the first kind for a general manifold $X$, 
and will return to
the second kind at the end when we discuss toric manifolds.
From now on, unless specified otherwise, $W_d$ will be
the first kind as defined above.

\item{(viii)}
We denote the equivariant hyperplane 
classes on $W_d$ by $\kappa_a$ (which are pullbacked
from the each of the $N_{d_a}$ to $W_d$).
We denote the equivariant hyperplane classes on $Y$ by $H_a$ 
(which are pullbacked from each of the $\P^{n_a}$ to $Y$).
We use the same notations for their restrictions to $X$.
We write $\kappa\cdot\zeta=\sum_a \kappa_a\zeta_a$,
$H\cdot t=\sum_a H_a t_a$, $d\cdot t=\sum_a d_a t_a$,
where the $t$ and $\zeta$ are formal variables.

\subsec{Localization on stable map moduli}

Clearly we have the commutative diagram:
\eqn\CommSquareII{
\matrix{
F_r & {\br i_r\over \longrightarrow} & M_d(X)\cr
e^Y\downarrow &  & \downarrow \varphi\cr
Y_r & {\br j_r\over \longrightarrow} & W_d.}
}
Let $\varphi:M_d(X)\ra W_d,~e^Y:F_r\ra Y_r$ play the respective roles
of $f:X\ra Y,~f':X'\ra Y'$ in functorial localization. 
Then it follows that

\lemma{Given a cohomology class $\omega$ on $M_d(X)$,
we have the following equality on $Y_r\cong Y$ for $0\preceq r\preceq d$:
$${j_r^*\varphi_*(\omega\cap LT_d(X))\over e_G(Y_r/W_d)}
=e^Y_*\left({i_r^*\omega\cap[F_r]^{vir}\over e_G(F_r/M_d(X))}\right).$$
}
\thmlab\SquareLemma

Following \LLYI, one can easily compute the
Euler classes $e_G(Y_r/W_d)$, and they are given 
as follows. For $d=(d_1,..,d_m)$, $r=(r_1,...,r_m)\preceq d$,
we have
$$e_G(Y_r/W_d)=\prod_{a=1}^m\prod_{i=0}^{n_a}{\prod_{k=0}^{d_a}}
_{k\neq r_a}
(H_a-\lambda_{a,i}-(k-r_a)\alpha)$$
where the $\lambda_{a,i}$ are the $T$ weights of $\P^{n_a}$.
Note that $e^Y$ is the composition of $e^X:F_r\ra X$
with $\tau:X\ra Y=Y_r$. Thus
$$e^Y_*=\tau_*e^X_*.$$
It follows that

\lemma{Given a cohomology class $\omega$ on $M_d(X)$,
we have the following equality on $Y_r$ for $0\preceq r\preceq d$:
$$
\tau^*\left({j_r^*\varphi_*(\omega\cap LT_d(X))\over e_G(Y_r/W_d)}\right)
=e_T(X/Y)\cap e^X_*\left({i_r^*\omega\cap[F_r]^{vir}\over e_G(F_r/M_d(X))}\right).
$$
}
\thmlab\SquareLemmaII

Now if $\psi$ is a cohomology class on $M_{0,0}(d,X)$,
then for $\omega=\pi^*\psi$, we get
$i_0^*\omega=i_0^*\pi^*\psi=\rho^*\psi$.
It follows that
%$$
%\frac{j_0^*\varphi_*(\pi^*\psi\cap LT_d(X))}{e_G(Y_0/W_d)}
%=e^Y_*\left(\frac{\rho^*\psi\cap LT_{0,1}(d,X)}{e_G(F_0/M_d(X))}\right).
%$$

\lemma{Given a cohomology class 
$\psi$ on $M_{0,0}(d,X)$,
we have the following equality on $X$:
$$
\tau^*\left({j_0^*\varphi_*(\pi^*\psi\cap LT_d(X))\over e_G(Y_0/W_d)}\right)
=e_T(X/Y)\cap 
e^X_*\left(\frac{\rho^*\psi\cap LT_{0,1}(d,X)}{e_G(F_0/M_d(X))}\right).
$$
}
\thmlab\SquareLemmaIII

\lemma{
For $r\neq 0, d$,
$$e_G(F_r/M_d(X))=
\alpha(\alpha+p_0^*c_1(L_r))\cdot
\alpha(\alpha-p_\infty^*c_1(L_{d-r})).
$$ 
For $r=0,d$, 
$$e_G(F_0/M_d(X))=\alpha(\alpha-c_1(L_d)), ~~~~
e_G(F_d/M_d(X))=\alpha(\alpha+c_1(L_d)).$$
}
\thmlab\NormalBundle
The computation 
done in section 2.3 of \LLYI~ and in section 3 of \LLYII~
(see also references there),
for the normal bundles $N_{F_r/M_d(X)}$,
makes no use of the convexity assumption on $TX$. Therefore it
carries over here with essentially no change.

%Combining the two preceding lemmas, we get
%the following equality on $Y_0$:
%$$\frac{j_0^*\varphi_*(\omega\cap LT_d(X))}{e_G(Y_0/W_d)}
%=e^Y_*\left(\frac{i_0^*\omega\cap [F_0]^{vir}}{\alpha(\alpha-c_1(L_d))}\right).
%$$

\subsec{From gluing identity to Euler data}

Fix a $T$-equivariant multiplicative class $b_T$.
Fix a $T$-equivariant bundle of the form $V=V^+\oplus V^-$,
where $V^\pm$ are respectively the convex/concave bundles on $X$.
We assume that
$$\Omega:={b_T(V^+)\over b_T(V^-)}$$
is a well-defined invertible class on $X$.
By convention, if $V=V^\pm$ is purely convex/concave,
then $\Omega=b_T(V^\pm)^{\pm1}$.
Recall that the bundle $V\ra X$ induces the bundles
$$V_d\ra M_{0,0}(d,X),~~
U_d\ra M_{0,1}(d,X),~~
\cU_d\ra M_d(X).$$
Moreover, they are related by $U_d=\rho^* V_d$, $\cU_d=\pi^* V_d$,
Define linear maps 
$$i_r^{vir}: A^G_*(M_d(X))\ra A^T_*(X)(\alpha),~~~~
i_r^{vir}\omega:=
e^X_*\left({i_r^*\omega\cap[F_r]^{vir}\over e_G(F_r/M_d(X))}\right).
$$

%Throughout this section, we denote
%$$
%Q:~~Q_d:=\varphi_*(\pi^* b_T(V_d)\cap LT_d(X)).
%$$
%If $\omega$ is a homology class on $W_d$,
%we write
%$$
%i_r^*\omega^v:=e_T(X/Y)^{-1}~\tau^*\left({j_r^*\omega\over
%e_G(Y_r/W_d)}\right),
%$$
%which is a class on $X$. By Lemma \SquareLemmaIII, we have
%$$
%i_r^*Q_d^v:=
%e^X_*\left(\frac{\rho^*\psi\cap LT_{0,1}(d,X)}{e_G(F_0/M_d(X))}\right).
%$$

\theorem{For $0\preceq r \preceq  d$, we have the following identity
in $A^T_*(X)(\alpha)$:
$$\Omega\cap i_r^{vir}\pi^*b_T(V_d)
=\overline{i_0^{vir}\pi^*b_T(V_r)}\cdot
i_0^{vir}\pi^*b_T(V_{d-r}).$$
}
\thmlab\GluingTheorem
\proof
For simplicity,
let's consider the case $V=V^+$.
The general case is entirely analogous.
The proof here is the one in \LLYI\LLYII,
but slightly modified to take into account the
new ingredient coming from the virtual class.

Recall that a point $(f,C)$ in $F_r\subset M_d$ comes
from gluing together
a pair of stable maps $(f_1,C_1,x_1),(f_2,C_2,x_2)$
with $f_1(x_1)=f_2(x_2)=p \in X$. From this, we get an exact sequence over $C$:
$$0\rightarrow f^*V\rightarrow f_1^*V\oplus f_2^*V\rightarrow
V|_p\rightarrow 0.$$ 
Passing to cohomology, we have   
$$0\rightarrow H^0(C, f^*V)\rightarrow H^0(C_1, f_1^*V)\oplus  H^0(C_2,
f_2^*V)
\rightarrow V|_p\rightarrow 0.$$
Hence we obtain an exact sequence of bundles
on $F_r$:
$$0\rightarrow i_r^*\cU_d\rightarrow U_r'\oplus U_{d-r}'\rightarrow
{e^X}^*V\rightarrow 0.$$
Here $i_r^*\cU_d$ is the restriction to $F_r$ of the bundle
 $\cU_d\ra M_d(X)$. And $U_r'$ is the pullback of the bundle
$U_r\ra M_{0,1}(d,X)$, and similarly for $U_{d-r}'$.
Taking the multiplicative class $b_T$, we get
the identity on $F_r$:
$${e^X}^*b_T(V)\cdot b_T(i^*_r\cU_d)=b_T(U_r')\cdot b_T(U_{d-r}').$$ 
We refer to this as
the {\it gluing identity}.

Now put
$$\omega={b_T(U_r)\over e_G(F_r/M_r(X))}\times
{b_T(U_{d-r})\over e_G(F_0/M_{d-r}(X))}\cap LT_{0,1}(r,X)\times LT_{0,1}(d-r,X)$$
From the commutative diagram \FiberSquare, we have the
identity:
$$
e^X_*\Delta^!(\omega)=
\Delta^*(e^X_r\times e^X_{d-r})_*(\omega).
$$
On the one hand is
$$\eqalign{
\Delta^*(e^X_r\times e^X_{d-r})_*(\omega)
&=(e^X_r)_*{b_T(U_r)\cap LT_{0,1}(r,X)\over e_G(F_r/M_r(X))}~\cdot~
(e^X_{d-r})_*{b_T(U_{d-r})\cap LT_{0,1}(d-r,X)\over e_G(F_0/M_{d-r}(X))}\cr
&=(e^X_r)_*{\rho^*b_T(V_r) \cap LT_{0,1}(r,X)\over e_G(F_r/M_r(X))}~\cdot~
(e^X_{d-r})_*{\rho^*b_T(V_{d-r}) \cap LT_{0,1}(d-r,X)\over e_G(F_0/M_{d-r}(X))}\cr
&=\overline{i_0^{vir}\pi^*b_T(V_r)}\cdot i_0^{vir}\pi^*b_T(V_{d-r}).
}
$$
On the other hand, applying the gluing identity, we have 
$$\eqalign{
e^X_*\Delta^!(\omega)
&=e^X_*\left({b_T(U_r')\over\alpha(\alpha+p_0^*c_1(L_r))}~\cdot~
{b_T(U_{d-r}')\over\alpha(\alpha-p_\infty^*c_1(L_{d-r}))}
\cap [F_r]^{vir}\right)\cr
&=e^X_*\left({ {e^X}^*b_T(V)\cdot 
i^*_rb_T(\cU_d)\cap [F_r]^{vir}\over e_G(F_r/M_d(X))}\right)\cr
&=b_T(V)\cap e^X_*\left({i^*_r b_T(\cU_d)\cap[F_r]^{vir}\over e_G(F_r/M_d(X))}\right)\cr
&=b_T(V)\cap i_r^{vir}\pi^*b_T(V_d).
}
$$
%the last equality being a consequence of Lemma \SquareLemma.
This proves our assertion.  $\Box$

Specializing the theorem to $b_T\equiv1$, we get

\corollary{
$i_r^{vir}1_d=\overline{i_0^{vir}1_r}\cdot i_0^{vir}1_{d-r}$
where $1_d$ is the identity class in on $M_d(X)$.
}

For a given convex/concave bundle $V$ on $X$,
and multiplicative  class $b_T$, we put
$$\eqalign{
A^{V,b_T}(t)=A(t)&:=e^{-H\cdot t/\alpha}\sum_d A_d~e^{d\cdot t}\cr
A_d&:=i_0^{vir}\pi^*b_T(V)=
e^X_*\left(\frac{\rho^*b_T(V_d)\cap LT_{0,1}(d,X)}{e_G(F_0/M_d(X))}\right).
}
$$
Here we will use the convention that $A_0=\Omega$,
and the sum is over all $d=(d_1,...,d_m)\in\Z_+^m$.
When the reference to $V,b_T$ is clear, we'll drop them
from the notations.
The special case in the corollary will play an important role.
So we introduce the notation:
$$
\One(t):=e^{-H\cdot t/\alpha}\sum_d \One_d e^{d\cdot t},~~~~
\One_d=i_0^{vir}1_d.
$$

By the preceding theorem and Lemma \SquareLemma, it follows 
immediately that for
$\omega=\varphi_*(\pi^* b_T(V_d)\cap LT_d(X))$, we have
$$\eqalign{
\int_{W_d} \omega\cap e^{\kappa\cdot\zeta}
&=\sum_{0\preceq r\preceq d}\int_{Y_r} 
{j_r^*\omega\over e_G(Y_r/W_d)}~
e^{(H+r\alpha)\cdot\zeta}\cr
&=\sum_r\int_{Y_r} 
\tau_* i_r^{vir}\pi^*b_T(V_d)~
e^{(H+r\alpha)\cdot\zeta}\cr
&=\sum_r\int_X
i_r^{vir}\pi^*b_T(V_d)~
e^{(H+r\alpha)\cdot\zeta}\cr
&=\sum_r\int_X
\Omega^{-1}\cap\bar A_r \cdot A_{d-r}~
e^{(H+r\alpha)\cdot\zeta}
~~~~~~(Theorem~\GluingTheorem).
}$$
Since $\omega\in A^G_*(W_d)$, hence
$\int_{W_d} \omega\cap c\in A_*^G(pt)=\C[\cT^*,\alpha]$
for all $c\in A_G^*(W_d)$, it follows that both sides of
the eqn. above lie in $\cR[[\zeta]]$.
This motivates the following (cf. \GiventalI)

\definition{Let $\Omega\in A_T^*(X)$, invertible.
We call a power series of the form
$$B(t):=e^{-H\cdot t/\alpha}\sum_d B_d ~e^{d\cdot t},~~~~
B_d\in A_*^T(X)(\alpha)
$$
an $\Omega$-Euler series if
$\sum_{0\preceq r\preceq d}\int_X
\Omega^{-1}\cap\bar B_r \cdot B_{d-r}~
e^{(H+r\alpha)\cdot\zeta}
\in\cR[[\zeta]]$ for all $d$.
}

Thus we have seen above that an elementary
consequence of the gluing identity in Theorem \GluingTheorem~
is that
\corollary{$A^{V,b_T}(t)
=e^{-H\cdot t/\alpha}\sum_d i_0^{vir}\pi^* b_T(V_d) ~e^{d\cdot t}$
is an Euler series.}

\definition{\LLYI~Let $\Lambda\in A_T^*(Y)$. We call a sequence
$P:~P_d\in A_G^*(W_d)$ an $\Lambda$-Euler data if
$$\Lambda\cdot j_r^*P_d=\overline{j_0^*P_r}\cdot j_0^*P_{d-r},
~~~0\preceq r\preceq d.$$
}
\thmlab\EulerData

Let $P$ be an $\Lambda$-Euler data such that $\tau^*\Lambda$ is
invertible. Then we have
\eqn\dumb{
\tau^*\Lambda\cdot\tau^*j_r^*P_d\cap i_r^{vir}1_d
=\tau^*\overline{j_0^*P_r}\cdot\tau^* j_0^*P_{d-r}
\cap \overline{i_0^{vir}1_r}\cdot i_0^{vir}1_{d-r}.
}
By Lemma \SquareLemmaII,
$$\eqalign{
\tau^*j_r^*P_d\cap i_r^{vir}1_d
&=\tau^*j_r^*P_d\cdot
e_T(X/Y)^{-1}\cap\tau^*\left({j_r^*\varphi_*LT_d(X)\over
e_G(Y_r/W_d)}\right)\cr
&=e_T(X/Y)^{-1}\cap\tau^*\left({j_r^*\varphi_*(\varphi^*P_d\cap LT_d(X))\over
e_G(Y_r/W_d)}\right)\cr
&=i_r^{vir}\varphi^* P_d.
}$$
Thus \dumb~ becomes
$$
\tau^*\Lambda\cap i_r^{vir}\varphi^*P_d
=\overline{i_0^{vir}\varphi^*P_r}\cdot i_0^{vir}\varphi^* P_{d-r}.
$$
(cf. Theorem \GluingTheorem.)
From this we get, as before,
$$
\int_{W_d} \varphi_*LT_d(X)\cap P_d\cdot e^{\kappa\cdot\zeta}
=\int_X\tau^*\Lambda^{-1}\cap 
\overline{i_0^{vir}\varphi^* P_r}\cdot
i_0^{vir}\varphi^* P_{d-r}~
e^{(H+r\alpha)\cdot\zeta}
\in\cR[[\zeta]].
$$
Therefore, that
$$
B(t)=
e^{-H\cdot t/\alpha}\sum_d i_0^{vir}\varphi^*P_d~e^{d\cdot t}
$$
is an Euler series, is just
an elementary consequence of
the Euler data identity.
More generally, we have

\theorem{Let $P$ be an $\Lambda$-Euler data as before, and let
$\O(t)=
e^{-H\cdot t/\alpha}\sum_d O_d~e^{d\cdot t}$
be any $\Omega$-Euler series. Then
$$B(t)=
e^{-H\cdot t/\alpha}\sum_d 
\tau^*j_0^*P_d\cap O_d
~e^{d\cdot t}$$
is an $\Omega\cdot \tau^*\Lambda$-Euler series.}
\thmlab\Multiplicative
\proof
Define $P'_d$ on $W_d$ by setting
$$j_r^*P'_d:=\tau_*(\Omega^{-1}\bar O_r\cdot O_{d-r})\cap e_G(Y_r/W_d).$$
By the localization theorem, this defines a class on $W_d$.
Moreover, we have
$$
\int_{W_d}P'_d\cap e^{\kappa\cdot\zeta}
=\sum_r\int_X\Omega^{-1} \bar O_r\cdot O_{d-r}~e^{(H+r\alpha)\cdot\zeta}
\in \cR.
$$
It follows that $P'_d\in A_*^G(W_d)\otimes\cR$
(see proof of Lemma 2.15 \LLYI).
Now
$$
\int_{W_d}P_d\cap P'_d
~e^{\kappa\cdot\zeta}
=\sum_r\int_X\Omega^{-1} \tau^*\Lambda^{-1}~
\overline{(\tau^*j_0^*P_r\cap O_r)}\cdot
(\tau^*j_0^*P_{d-r}\cap O_{d-r})~e^{(H+r\alpha)\cdot\zeta},
$$
which lies in $\cR$ because
$P_d\cap P'_d$ lies in $A_*^G(W_d)\otimes\cR$.
$\Box$

Note that if $O_d=\One_d$, then $P'_d$ in
the proof above is just $\varphi_*LT_d(X)$.
For explicit examples of Euler data, see \LLYI\LLYII.

\subsec{From Euler data to intersection numbers}

Again, fix the data $V,b_T$ as before.
From now on we write $e^X$ simply as $e$.
We recall the notations
$$\eqalign{
A^{V,b_T}(t)=A(t)&=e^{-H\cdot t/\alpha}
\sum_d A_d~e^{d\cdot t},\cr
A_d&=i^{vir}_0\pi^*b_T(V_d)=
e_*^X\left({\rho^* b_T(V_d)\cap LT_{0,1}(d,X)\over e_G(F_0/M_d(X))}\right)\cr
K_d^{V,b}=K_d&=\int_{LT_{0,0}(d,X)} b(V_d)\cr
\Phi^{V,b}=\Phi&=\sum K_d~e^{d\cdot t}.
}
$$

\theorem{(i) $deg_\alpha A_d\leq-2$.
\item{(ii)} If for each $d$ the class $b_T(V_d)$ has homogeneous
degree the same as the degree of
$M_{0,0}(d,X)$, then in
 the nonequivariant limit we have
$$\eqalign{
\int_X e^{-H\cdot t/\alpha}
A_d&=\alpha^{-3}(2-d\cdot t)K_d\cr
\int_X\left(A(t)-e^{-H\cdot t/\alpha}\Omega\right)&
=\alpha^{-3}(2\Phi-\sum t_i {\partial\Phi\over\partial t_i}).
}$$
}
\thmlab\KdTheorem
\proof
By definition,
$$
A_d=e_*\left({\rho^* b(V_d)\cap LT_{0,1}(d,X)\over
e_G(F_0/M_d(X))}\right).
$$ 
So assertion (i) follows immediately from this formula 
Lemma \NormalBundle.

The second equality in assertion (ii) follows from
the first equality in (ii). Now consider
$$\eqalign{
I&:=\int_X e^{-H\cdot t/\alpha} A_d\cr
&=\int_{LT_{0,1}(d,X)} e^{-e^*H\cdot t/\alpha}
{\rho^* b(V_d)\over
e_{\C^\times}(F_0/M_d(X))}\cr
&=\int_{LT_{0,0}(d,X)}  b(V_d)~
\rho_*\left({e^{-e^*H\cdot t/\alpha}\over
e_{\C^\times}(F_0/M_d(X))}\right).
}$$
Now $b(V_d)$ has homogeneous degree the same as the dimension
of $LT_{0,0}(d,X)$.
The second factor in the last integrand contributes
a scalar factor given by integration
over a fiber $E$
of $\rho$. 
By Lemma \NormalBundle,
the degree 1 term in the second factor is
${-e^*H\cdot t\over \alpha^3}+{c\over \alpha^3}$
where $c=c_1(L_d)$.

Now the line bundle $L_d$ on $M_{0,1}(d,X)$
is the restriction of the universal bundle $L_d'$
on $M_{0,1}(d,Y)$,
and the map $\rho:M_{0,1}(d,X)\ra M_{0,0}(d,X)$,
is the restriction of the forgetting map
$\rho':M_{0,1}(d,Y)\ra M_{0,0}(d,Y)$. For the latter,
we can choose a smooth fiber $E'\cong\P^1$ so that
$$\int_{E'} i_{E'}^*c_1(L_d')=\int_{E'} c_1(TE')= 2.$$
Since $\rho'$ is flat, 
$$\int_E i_E^*c_1(L_d)=\int_{E'}i_{E'}^*c_1(L_d')=2.$$
Restricting to a fiber $E$
say over $(C,f)\in M_{0,0}(d,X)$, the
evaluation map $e$ is equal to $f$, which is
a degree $d$ map $E\ra X$.
It follows that
$$\int_E e^*H=d.$$
So we have
$$I=(-{d\cdot t\over \alpha^3}+{2\over \alpha^3})K_d.~~~~\Box$$

\theorem{
More generally suppose $b_T$ is an equivariant
multiplicative class of the form
$$b_T(V)=x^r+x^{r-1}b_1(V)+\cdots+b_r(V),~~rk~V=r$$
where $x$ is a formal variable, $b_i$ is a 
class of degree $i$.
Suppose $s:=rk~V_d-exp.~dim~
M_{0,0}(d,X)\geq0$ is independent of $d\succ0$. Then
in the nonequivariant limit,
$$\eqalign{
{1\over s!}\left({d\over dx}\right)^s|_{x=0}\int_X e^{-H\cdot t/\alpha}
A_d&=\alpha^{-3}x^{-s}(2-d\cdot t)K_d\cr
{1\over s!}\left({d\over dx}\right)^s|_{x=0}
\int_X\left(A(t)-e^{-H\cdot t/\alpha}\Omega\right)&
=\alpha^{-3}x^{-s}(2\Phi-\sum t_i {\partial\Phi\over\partial t_i}).
}$$
}
\thmlab\KdTheoremII
\proof
The proof is entirely analogous to (ii) above. $\Box$

In the case of $b_T(V)=1$,
one can improve the $\alpha$ degree estimates for $A_d=\One_d$
given by Theorem \KdTheorem~(i). 

\lemma{
For all $d$,
$$deg_\alpha~ \One_d\leq min(-2,-\bra c_1(X),d\ket).$$
}
\thmlab\DegreeBound
\proof
If $\bra c_1(X),d\ket\leq 2$, then the assertion 
is a special case Theorem \KdTheorem~(i). So suppose 
that $\bra c_1(X),d\ket> 2$. 
The class $LT_{0,1}(d,X)$ is of dimension
$s=exp.dim~M_{0,1}(d,X)=\bra c_1(X),d\ket+dim~X-2$. 
Let $c=c_1(L_d)$. Then $c^k\cap LT_{0,1}(d,X)$ is
of dimension $s-k$, and so $e_*(c^k\cap LT_{0,1}(d,X))$
lies in the group $A^T_{s-k}(X)$. But this group is zero
unless $s-k\leq dim~X$ or $k\geq s-dim~X=\bra c_1(X),d\ket-2$.
Now by Lemma \NormalBundle, it follows that
$$
\One_d=i_0^{vir}1_d=\sum_{k\geq\bra c_1(X),d\ket -2}
{1\over \alpha^{k+2}}~e_*(c^k\cap LT_{0,1}(d,X)).
$$
This completes the proof. $\Box$

\remark{The entire theory discussed in this section
obviously specializes to the case $T=1$,
hence applies to
any projective manifold $X$.
}

\newsec{Linking}

\definition{A projective $T$-manifold $X$ is called a 
balloon manifold if $X^T$ is finite, and if
for $p\in X^T$, the weights of
the isotropic representation $T_pX$ are pairwise
linearly independent.}

The second condition in the definition is known as
the GKM condition \GKM. We will assume that our balloon
manifold has the property that if $p,q\in X^T$
such that $c(p)=c(q)$ for all $c\in A^1_T(X)$, then $p=q$.
{\it From now on, unless stated otherwise,
$X$ will be a balloon manifold
with this property. }
If two fixed points $p,q$ in $X$ are connected by
a $T$-invariant 2-sphere, then we call that 2-sphere a balloon
and denote it by $pq$.
For examples and the basic facts we need to use about
these manifolds, see \LLYII~ and references there.
All the results in sections 5-6 in \LLYII~ are
proved for balloon manifolds without any convexity
assumption, and are therefore also applicable here.
We will quote the ones we need here without proof,
but with only slight change in notations and terminology.

\definition{Two Euler series $A,B$ are linked if
for every balloon $pq$ in $X$ and
every $d=\delta[pq]\succ0$, the function
$(A_d-B_d)|_p\in \C(\cT^*)(\alpha)$
is regular
at $\alpha={\lambda\over\delta}$
where $\lambda$ is the 
weight on the tangent line $T_p(pq)\subset T_pX$.
}

\theorem{(Theorem 5.4 \LLYII) 
Suppose $A,B$ are linked Euler series satisfying
the following properties: for $d\succ0$,
\item{(i)} For $p\in X^T$,
every possible pole of
$(A_d-B_d)|_p$
is a scalar multiple of a weight on $T_pX$.
\item{(ii)} 
$deg_\alpha (A_d-B_d)\leq -2$.
\item{} Then we have $A=B$.}
\thmlab\Uniqueness

%In our applications later, the situation is
%better then the conditions (i)-(ii) demand.
%We will have two Euler series $A,B$ such that
%$A_d,B_d$ separately, rather than just $A_d-B_d$,
%will satisfy both conditions (i)-(ii) at the
%outset. In this situation, to prove that $A=B$,
%it suffices to prove that they are linked.

\theorem{(Theorem 6.6 \LLYII)
Suppose that $A,B$ are two linked Euler series having property (i) of
the preceding theorem.
Suppose that $deg_\alpha A_d\leq-2$
for all $d\succ0$,
and that there exists
power series $f\in\cR[[e^{t_1},..,e^{t_m}]]$, 
$g=(g_1,..,g_m)$, $g_j\in\cR[[e^{t_1},..,e^{t_m}]]$,
without constant terms, such that
\eqn\Asymptotic{e^{f/\alpha}B(t)=\Omega- 
\Omega {H\cdot (t+g)\over\alpha}
+O(\alpha^{-2})}
when expanded in powers of $\alpha^{-1}$.  Then 
$$A(t+g)=e^{f/\alpha}~B(t).$$
}
\thmlab\Applications

The change of variables effected by $f,g$ above is
an abstraction of what's known as mirror transformations \CDGP.

\theorem{Let $p\in X^T$, 
$\omega\in A^*_T(M_{0,1}(d,X))[\alpha]$, and consider
$i_p^* e_*\left(
{\omega\cap LT_{0,1}(d,X)\over e_G(F_0/M_d(X))}\right)\in\C(\cT^*)(\alpha)$
as a function of $\alpha$. Then
\item{(i)} Every possible pole of the function is a scalar multiple of
a weight on $T_pX$.
\item{(ii)} Let $pq$ be a balloon in $X$,
and $\lambda$ be the weight on the tangent line $T_p(pq)$.
If $d=\delta[pq]\succ0$,
then the pole of the function at $\alpha=\lambda/\delta$ 
is of the form
$$
e_T(p/X)~{1\over\delta}{1\over\alpha(\alpha-\lambda/\delta)}
{i_F^*\omega\over e_T(F/M_{0,1}(d,X))}
$$
where $F$ is the (isolated) fixed point $(\P^1,f_\delta,0)\in M_{0,1}(d,X)$
with $f_\delta(0)=p$, and $f_\delta:\P^1\ra X$ maps by a $\delta$-fold cover
of $pq$.
}
\thmlab\SimplePole
\proof
Consider the commutative diagram
$$
\matrix{
\{F\} & {\br i_F\over\longrightarrow}& 
M_{0,1}(d,X)\cr
e'\downarrow &  &~~e\downarrow\cr
p & {\br i_p\over\longrightarrow} & X
}
$$
where $e$ is the evaluation map,
$\{F\}$ are the fixed point components in $e^{-1}(p)$,
$e'$ is the restriction of $e$ to $\{F\}$, and $i_F,i_p$
are the usual inclusions. By functorial localization
we have, for any $\beta\in A^*_T(M_{0,1}(d,X))(\alpha)$,
\eqn\dumb{\eqalign{
i_p^* e_*(\beta\cap LT_{0,1}(d,X))
&=e_T(p/X)
~\sum_F e'_*\left(
{i_F^*\beta\cap[F]^{vir}\over e_T(F/M_{0,1}(d,X))}\right)\cr
&=e_T(p/X)~
\sum_F\int_{[F]^{vir}}
{i_F^*\beta\over e_T(F/M_{0,1}(d,X))}.
}}
We apply this to the class
$$
\beta=
{\omega\over e_G(F_0/M_d(X))}
={\omega\over \alpha(\alpha-c)}
$$
where $c=c_1(L_d)$ (cf. Lemma \NormalBundle).
For (i), we will show that a pole of the sum \dumb~
is at either $\alpha=0$ or $\alpha=\lambda'/\delta'$
for some tangent weight $\lambda'$ on $T_pX$.
For (ii), we will show that only one $F$ in the sum \dumb~
contributes to the pole at $\alpha=\lambda/\delta$,
that the contributing $F$ is the isolated fixed point
$(\P^1,f_\delta,0)$ as asserted in (ii), and that the contribution has
the desired form.

A fixed point $(C,f,x)$ in $e^{-1}(p)$ is such that 
$f(x)=p$, and that the image curve $f(C)$ 
lies in the union of the $T$-invariant
balloons in $X$.
The restriction of the first Chern class $c$ to
an $F$ must be of the form
$$i_F^*c=c_F+w_F$$
where $c_F\in A^1(F)$, and $w_F\in\cT^*$
is the weight of the representation on the line $T_xC$
induced by the linear map $df_x:T_xC\ra T_pX$ (cf. \Kontsevich).
The image of $df_x$ is either 0 or a tangent line $T_p(pr)$
of a balloon $pr$.
Thus $w_F$ is either zero (when the branch $C_1\subset C$ containing $x$
is contracted), or $w_F=\lambda'/\delta'$
(when $C_1{\br f\over\ra} X$ maps by a $\delta'$-fold cover
of a balloon $pr$ with tangent weight $\lambda'$).
The class $e_T(F/M_{0,1}(d,X))$
is obviously independent of $\alpha$.
Since $c_F$ is nilpotent, a pole of the sum \dumb~
is either at $\alpha=0$ or $\alpha=w_F$ for some $F$.
This proves (i).

Now, an $F$ in the
sum \dumb~ contributes to the pole at $\alpha=\lambda/\delta$
only if $w_F=\lambda/\delta$.
Since the weights on $T_pX$ are pairwise linearly independent,
that $\lambda/\delta=\lambda'/\delta'$ implies that 
$\lambda=\lambda'$ and $\delta=\delta'$.
Since $d=\delta[pq]$, it
follows that the only fixed point contributing
to the pole at $\alpha=\lambda/\delta$ is $(C,f,x)$ where 
$C{\br f\over\ra} X$ maps by a $\delta$-fold cover
of the balloon $pq$ with $C\cong\P^1$ and $f(x)=0$.
This is an isolated fixed point, which we 
denote by $F=(\P^1,f_\delta,0)$.
It contributes to the sum \dumb~ the term
$$
\int_F{i_F^*\beta\over e_T(F/M_{0,1}(d,X))}
={1\over\delta}
{i_F^*\omega\over \alpha(\alpha-\lambda/\delta)}
{1\over e_T(F/M_{0,1}(d,X))}.
$$
Here $F$ is an orbifold point of order $\delta$,
and hence the integration contributes the factor $1/\delta$.
This proves (ii). $\Box$

%\corollary{In the same notations as in the preceding theorem (ii),
%if $i_F^*\omega$ vanishes at $\alpha=\lambda/\delta$
%then the  function
%$i_p^* e_*\left(
%{\omega\cap LT_{0,1}(d,X)\over e_G(F_0/M_d(X))}\right)$
%is regular there.
%}
%
\def\One{{\bf 1}}

Fix the data $V,b_T$ and a $\Lambda$-Euler data $P:~P_d$
such that
$$\tau^*\Lambda=\Omega:=b_T(V^+)/b_T(V^-).$$
We now discuss the interplay between four Euler series:
$A^{V,b_T}(t)$, $\One(t)$, and 
two others
$$\eqalign{
\O(t)&:=e^{-H\cdot t/\alpha}\sum
O_d~e^{d\cdot t}\cr
B(t)&:=e^{-H\cdot t/\alpha}\sum_d 
\tau^*j_0^*P_d\cap O_d
~e^{d\cdot t}
}$$
{\it where $\O(t)$ denote some unspecified Euler series linked
to $\One(t)$.} (In particular $\O(t)$ may be
specialized to $\One(t)$ itself.) That $B(t)$ is
an Euler series follows from Theorem \Multiplicative.

\corollary{Suppose that
at $\alpha=\lambda/\delta$ and $F=(\P^1,f_\delta,0)$,
we have $i_p^*j_0^*P_d=i_F^*\rho^* b_T(V_d)$
for all $d=\delta[pq]$.
Then $B(t)$ 
is linked to $A^{V,b_T}(t)$.
}
\proof
Since $\O(t)$ is linked to $\One(t)$ by assumption,
it follows trivially that
$$\eqalign{
B(t)&=e^{-H\cdot t/\alpha}\sum_d 
\tau^*j_0^*P_d\cap O_d
~e^{d\cdot t}\cr
C(t)&=e^{-H\cdot t/\alpha}\sum_d 
\tau^*j_0^*P_d\cap \One_d
~e^{d\cdot t}
}
$$
are linked. So it suffices to show that
$A(t)$ and $C(t)$ are linked.
Denote their
respective Fourier coefficients by $A_d$, $C_d$.
Then
\eqn\dumb{
i_p^*C_d - i_p^*A_d
=i_p^*j_0^*P_d\cdot 
i_p^* e_*\left(
{LT_{0,1}(d,X)\over e_G(F_0/M_d(X))}\right)
-i_p^* e_*\left(
{\rho^*b_T(V_d)\cap LT_{0,1}(d,X)\over e_G(F_0/M_d(X))}\right).
}
By Theorem \SimplePole~ (ii), this difference is regular
because the zero of the function
$i_p^*j_0^*P_d-i_F^*\rho^* b_T(V_d)$
cancels the {\it simple} pole
of each term in \dumb~ 
at $\alpha=\lambda/\delta$. $\Box$

We now formulate one of the main theorems of this paper.
It'll also give a more directly applicable form of
Theorem \Applications. Given the data
$V,b_T,\O(t),P$, and 
$$
B(t):=e^{-H\cdot t/\alpha}\sum_d 
\tau^*j_0^*P_d\cap O_d~ e^{d\cdot t},
$$
assume that the preceding corollary holds.
Suppose in addition, that
\item{(*)} For each $d$, we have the form
$$
\tau^*j_0^* P_d=\Omega\alpha^{\bra c_1(X),d\ket}
\left(a+(a'+ a''\cdot H)\alpha^{-1}+\cdots\right),
$$
for some $a,a',a_i''\in\C(\cT^*)$ (depending on $d$).
\item{(**)} For each $d$, we have the form
(written in cohomology $A^*_T(X)$):
$$
O_d=\alpha^{-\bra c_1(X),d\ket}
\left(b+(b'+ b''\cdot H)\alpha^{-1}+\cdots\right),
$$
for some $b,b',b_i''\in\C(\cT^*)$ (depending on $d$).

\theorem{Suppose that $A^{V,b_T}(t),B(t)$
are as in the preceding corollary. Under the assumptions (*)-(**),
there exist power series 
$f\in\cR[[e^{t_1},..,e^{t_m}]]$, 
$g=(g_1,..,g_m)$, $g_j\in\cR[[e^{t_1},..,e^{t_m}]]$,
without constant terms, such that
$$A^{V.b_T}(t+g)=e^{f/\alpha} B(t).$$
}
\thmlab\ApplicationsII
\proof
Recall that
$$
B(t):=e^{-H\cdot t/\alpha}\sum_d 
\tau^*j_0^*P_d\cap O_d
~e^{d\cdot t}.
$$
By the preceding corollary, $B(t)$ is linked to $A(t)$.
We will use the asymptotic forms (*)-(**) to explicitly
construct $f,g$ satisfying the condition \Asymptotic.
Our assertion then follows from Theorem \Applications.

By (*)-(**), the Fourier coefficient $B_d$, $d\succ0$, of $B(t)$ has
the form
$$
B_d=
\Omega\left(ab+(ab'+a'b)\alpha^{-1}+ (ab''+a''b)\cdot H\alpha^{-1}+\cdots\right)
$$
(and $B_0=\Omega$). 
Multiplying this by $e^{-H\cdot t/\alpha}=1-H\cdot t\alpha^{-1}+\cdots$,
and $e^{d\cdot t}$, and then sum over $d$,
we get the form
$$
B(t)=
\Omega\left(C+ (C'+C''\cdot H-C~H\cdot t)\alpha^{-1}+\cdots\right)
$$
where $C,C',C''_i\in\C(\cT^*)[[e^{t_1},..,e^{t_m}]]$
having constant terms $1,0,0$ respectively.
It follows that
$$
{e^{-C'/C\alpha}\over C} B(t)
=\Omega\left(1-(t-{C''\over C})\cdot H\alpha^{-1}+\cdots\right)
$$
So putting 
$f=-\alpha log~C-{C'\over C}$ and $g=-{C''\over C}$
yields the eqn. \Asymptotic. This completes the proof.  $\Box$

\corollary{The preceding theorem holds if we specialize
the choice of $O(t)$ to $\One(t)$, ie.
$$
B(t)=e^{-H\cdot t/\alpha}\sum_d 
\tau^*j_0^*P_d\cap \One_d
~e^{d\cdot t}.
$$
}
\proof
The preceding theorem holds for any
Euler series $\O(t)$ 
satisfying the condition (**)
linked to $\One(t)$.
Now by Lemma
\DegreeBound, $\One(t)$ satisfies condition (**);
and obviously it is also linked to itself.
$\Box$

\subsec{Linking values}

In this subsection, we continue using the notations
$V,b_T,\One(t),O(t),A(t)$, introduced above,
where $O(t)$ is linked to $\One(t)$.
We will apply Theorem \ApplicationsII~ to
the case when $b_T$ is the Euler class or
the Chern polynomial.

For simplicity, we will assume that $V$
has the following property: that there exist
nontrivial $T$-equivariant line bundles
$L_1^+,..,L_{N^+}^+;L^-_1,..,L^-_{N^-}$ on $X$ with
$c_1(L_i^+)\geq0$ and $c_1(L_j^-)<0$, such that
for any balloon $pq\cong\P^1$ in $X$ we have
$$V^\pm|_{pq}=\oplus_{i=1}^{N^\pm}L_i^\pm|_{pq}.$$
Note that $N^\pm=rk~V^\pm$.
We also require that 
\eqn\FactorizedForm{
\Omega:=b_T(V^+)/b_T(V^-)=\prod_i b_T(L^+_i)/\prod_j b_T(L^-_j).
}
In this case we call the list 
$(L_1^+,..,L_{N^+}^+;L^-_1,..,L^-_{N^-})$ the
splitting type of $V$. Note that $V$ is not assumed to
split over $X$.

\theorem{Let $b_T=e_T$ be the equivariant Euler class.
Let $pq$ be a balloon, $d=\delta[pq]\succ0$,
and $\lambda$ be the weight on the tangent line $T_q(pq)$. 
Let $F=(\P^1,f_\delta,0)$ be the fixed point, as in Theorem \SimplePole (ii).
Then
\eqn\dumb{i_F^*\rho^*b_T(V_d)=
\prod_i\prod_{k=0}^{\bra c_1(L^+_i),d\ket}
\left(c_1(L^+_i)|_p-k\lambda/\delta\right)\times
\prod_j\prod_{k=1}^{-\bra c_1(L^-_j),d\ket-1}
\left(c_1(L^-_j)|_p+k\lambda/\delta\right).
}
In particular, $A^{V,e_T}(t)$ 
is linked to the Euler series
$B(t)=e^{-H\cdot t/\alpha}\sum_d 
B_d~e^{d\cdot t}$ where
$$B_d=
O_d\cap
\prod_i\prod_{k=0}^{\bra c_1(L^+_i),d\ket}
(c_1(L^+_i)-k\alpha)\times
\prod_j\prod_{k=1}^{-\bra c_1(L^-_j),d\ket-1}
(c_1(L^-_j)+k\alpha).
$$
}
\thmlab\EulerClass
\proof 
Define 
$P:~~P_d\in A^*_G(W_d)$ by
$$
P_d:=\prod_i\prod_{k=0}^{\bra c_1(L^+_i),d\ket}
(\hat L^+_i-k\alpha)\times
\prod_j\prod_{k=1}^{-\bra c_1(L^-_j),d\ket-1}
(\hat L^-_j+k\alpha),
$$
where $\hat L^\pm_i\in A_G^*(W_d)$ denotes the canonical lifting
of $c_1(L^\pm_i)\in A_T^*(Y)$. Then $P$ is an $\Omega$-Euler data
(see section 2.2 \LLYI).
By Theorem \Multiplicative, it follows that
$B_d=\tau^* j_0^*P_d\cap O_d$ is
an Euler series. By (corollary to) Theorem \SimplePole,
$A(t)$ is linked to $B(t)$, provided that
eqn. \dumb~ holds.
We now prove eqn. \dumb.

We first consider a single convex line bundle $V=L$.
As before, the fixed point $F=(\P^1,f_\delta,0)$ in $M_{0,1}(d,X)$ is a
$\delta$-fold cover of the balloon $pq\cong\P^1$.
We can write it as
$$
f_\delta:\P^1\ra pq\cong\P^1,~~~
[w_0, w_1]\mapsto[w_0^\delta, w_1^\delta].
$$ 
Note that the $T$-action on $X$ induces the
standard rotation on $pq\cong \P^1$ with weight $\lambda$.
Clearly, we have 
$$
i_F^*\rho^*e_T(V_d)=i_{\rho(F)}^*e_T(V_d)=
e_T(i^*_{\rho(F)}V_d).
$$
The right hand side is the product the weights of the $T$ representation
on the vector space
$$
i^*_{\rho(F)}V_d=H^0(\P^1,f^*_\delta L)=H^0(\P^1,f^*_\delta\cO(l))
$$
where $l=\bra c_1(L),[pq]\ket$. 
Thus we get
(cf. section 2.4 \LLYI)
$$
e_T(i^*_{\rho(F)}V_d)
=\prod_{k=0}^{l\delta}(c_1(L)|_p-k\frac{\lambda}{\delta}).
$$
This proves \dumb~ for a single convex line bundle.

Similarly for a concave line bundle $V=L$, if its restriction to the
balloon $pq$ is $\cO(-l)$ with $-l=\bra c_1(L), [pq]\ket$, then 
$$
e_T(i^*_{\rho(F)}V_d)
=\prod_{k=1}^{l\delta-1}(c_1(L)|_p+k\frac{\lambda}{\delta}).
$$
This is \dumb~ for a single concave line bundle.
The general case can clearly be obtained by taking products. 
$\Box$

A parallel argument for $b_T=$ the Chern polynomial yields

\theorem{Let $b_T=c_T$ be the equivariant Chern polynomial,
with the rest of the notations as in the preceding theorem.
Then
$$
i_F^*\rho^* c_T(V)
=\prod_i\prod_{k=0}^{\bra c_1(L^+_i),d\ket}
\left(x+c_1(L^+_i)|_p-k\lambda/\delta\right)\times
\prod_j\prod_{k=1}^{-\bra c_1(L^-_j),d\ket-1}
\left(x+c_1(L^-_j)|_p+k\lambda/\delta\right).
$$
In particular, $A^{V,e_T}(t)$ 
is linked to the Euler series
$B(t)=e^{-H\cdot t/\alpha}\sum_d 
B_d~e^{d\cdot t}$ where
$$B_d=
O_d\cap
\prod_i\prod_{k=0}^{\bra c_1(L^+_i),d\ket}
(x+c_1(L^+_i)-k\alpha)\times
\prod_j\prod_{k=1}^{-\bra c_1(L^-_j),d\ket-1}
(x-c_1(L^-_j)+k\alpha).
$$
}
\thmlab\ChernPolynomial

By Theorem \ApplicationsII, we can therefore
compute $A(t)=A^{V,b_T}(t)$ in terms of the Euler series $B(t)$
given above, provided that the Euler data $P$ and
the Euler series $O(t)$ both have
the appropriate asymptotic forms (*)-(**) required by
Theorem \ApplicationsII. 

\corollary{
Let $b_T$ be either $e_T$ or $c_T$. Suppose that
\eqn\ChernClassBound{
c_1(V^+)-c_1(V^-)\leq c_1(X).
}
Then the condition (*) holds for the Euler data $P$
in the two preceding theorems. In this case, if $O(t)$ is any
Euler series linked to $\One(t)$ and satisfies
condition (**), then Theorem \ApplicationsII~
applies to compute $A^{V,b_T}(t)$ in terms of $O(t)$ and $P$.}
\thmlab\ApplicationsIII
\proof
The Euler data $P$ in either of the preceding theorems
has the form: for each $d\succ0$,
$$
\tau^*j_0^* P_d=\Omega\alpha^{\bra c_1(V^+)-c_1(V^-),d\ket-N_-}
\left(a+(a'+ a''\cdot H)\alpha^{-1}+\cdots\right),
$$
for some $a,a',a_i''\in\C(cT^*)$ (depending on $d$).
By assumption,
$$\bra c_1(V^+)-c_1(V^-),d\ket\leq \bra c_1(X),d\ket.$$
This implies that $P$ satisfies the condition (*).
$\Box$

This result shows that if $\Omega=b_T(V^+)/b_T(V^-)$
has a certain factorized form \FactorizedForm, 
and if there is a suitable bound \ChernClassBound~
on first Chern classes, then $A(t)=A^{V,b_T}(t)$
is computable in terms of the $\One(t)$
(or a suitable Euler series $O(t)$ linked to it).
Note that even though $\One(t)$ is not known explicitly
in closed form in general, it is {\it universal} in the sense
that it is natural and is independent of any choice of $V$ or $b_T$.
Its Fourier coefficients also happen to be related to
the universal line bundle on $M_{0,1}(d,X)$.
In the next section, 
we specialize
$O(t)$ to something quite explicit.
We also discuss some other ways
to compute $A(t)$. 
We consider situations in which
the first Chern class bound and
the factorization condition on $\Omega$ can be removed.

\newsec{Applications and Generalizations}

Throughout this section, we continue to use the
same notations: $V,b_T,\Omega,A(t),...$.

\subsec{Inverting $\One_d$}

Suppose
$\One_d$ is invertible for all $d$.
Then obviously, there exist unique $B_d\in A^*_T(X)(\alpha)$
such that
$$A(t)=e^{-H\cdot t/\alpha}\sum B_d\cap  \One_d~e^{d\cdot t}.$$
In particular this says that
for  $d=\delta[pq]$, $F=(\P^1,f_\delta,0)$, we must have
\eqn\dumb{
i_p^*B_d=i_F^*\rho^*b_T(V_d)}
at $\alpha=\lambda/\delta$.
By Theorem \Uniqueness, the $B_d$ are the unique classes in
$A^*_T(X)(\alpha)$ such that
\item{(i)} eqn. \dumb~ holds.
\item{(ii)} $deg_\alpha~B_d\cap \One_d\leq -2$.
\item{(iii)} 
$e^{-H\cdot t/\alpha}\sum B_d\cap  \One_d~e^{d\cdot t}$ is
an $\Omega$-Euler series.

In other words these algebraic conditions completely
determine the $B_d$. Thus in principle the $B_d$
can be computed in terms of the classes
$\One_d$ and the linking
values \dumb. The point is that {\it this is true whether or not
the bound \ChernClassBound~ or the factorized form
$\Omega$ \FactorizedForm~ holds.}
Here are a few examples. 

\bigskip
\noindent {\it Example 1.} $X=Y$ is a product of projective space
with the maximal torus action.
In this case, 
$$\One_d={1\over e_G(Y_0/W_d)}$$
which is given explicitly in section 2.
We also have $B_d=j_0^*\varphi_*(\pi^*b_T(V_d)\cap LT_d(X))\in A^*_T(X)[\alpha]$
(cf. Lemma \SquareLemmaII). Finding the $B_d$ explicitly
amounts to finding polynomials in $H_a,\alpha$
with the prescribed values \dumb, and the degree bound (ii).
This is a linear problem!
This approach is particularly useful for computing
$b_T(V_d)$ for nonsplit bundles $V$ (e.g. $V=TX$),
or for bundles where the bound \ChernClassBound~ fails
(e.g. $\cO(k)$ on $\P^n$ with $k>n+1$).

\bigskip
\noindent {\it Example 2.} Suppose $X$ is a balloon manifold
such that every balloon $pq$ generates the integral
classes in $A_1(X)$. Then every
integral class $d\in A_1(X)$ is of the form $\delta[pq]$
(e.g. Grassmannians).
We claim that, in this case, $\One_d$ is invertible
for all $d$.
It suffices to show that
$i_p^* \One_d$ is nonzero for every fixed point $p$
in $X$. Given $p$, we know that there are $n=dim~X$
other fixed points $q$ joint to $p$ by balloons $pq$.
Pick such a $q$. Then
$d=\delta[pq]$ for some $\delta$. It follows from
Theorem \SimplePole~ that the function
$i_p^* \One_d$ has a {\it nontrivial} simple pole
at $\alpha=\lambda/\delta$ where $\lambda$ is
the weight on the tangent line $T_p(pq)$.
This completes the proof.

Obviously, we can take product of these examples
and still get invertible $\One_d$ for
the product manifold.

\subsec{Toric manifolds}

Let $X$ be a toric manifold of dimension $n$. Denote by $D_a$, $a=1,..,N$,
the $T$-invariant divisors in $X$. We denote by the same notations the
equivariant homology classes they represent. Recall that 
\Audin\Cox\Musson~
$X$ can be represented as an orbit space
$$X=(\Gamma-Z)/K$$
where $K$ is an algebraic torus of dimension $N-n$,
$\Gamma=\C^N$ is a linear representation of $K$, 
and $Z$ is a $K$-invariant
monomial variety
of $\C^N$,
all determined by the fan of $X$. The $T$ action
on the orbit space is induced by  $(\C^\times)^N$
acting on $\Gamma$ by the usual scaling.
Define
\eqn\FormulaO{
O(t)= e^{-H\cdot t/\alpha}\sum_d O_d~e^{d\cdot t},~~~~
O_d:={\prod_{\bra D_a,d\ket<0}
\prod_{k=0}^{-\bra D_a,d\ket-1}(D_a+k\alpha)\over
\prod_{\bra D_a,d\ket\geq0}
\prod_{k=1}^{\bra D_a,d\ket}(D_a-k\alpha)}.
}
We will prove that $O(t)$ is a $1$-Euler series.

First we recall a construction in \MP\Witten. Given
an integral class $d\in A_1(X)$,
let $\Gamma_d=\oplus_a H^0(\P^1,\cO(\bra D_a,d\ket))$. Let
$K$ act on $\Gamma_d$ by $\phi_a\mapsto t^{\lambda_a}\phi_a$
where the $\lambda_a$ are the same weights with which $K$
acts on $\Gamma$. Let 
$$
Z_d=\{\phi\in\Gamma_d|
\phi(z,w)\in Z,~\forall(z,w)\in\C^2\}.
$$
(Note that $\phi$ here is viewed as a polynomial map $\C^2\ra\C^N$.)
It is obvious that $Z_d$ is $K$-invariant.
Define the orbit space
$$\cW_d:=(\Gamma_d-Z_d)/K.$$

\item{(i)} If not empty, $\cW_d$ is a toric manifold of dimension
$$dim~\cW_d=\sum_a'(\bra D_a,d\ket+1)-dim~K$$
where $\sum'_a$ means summing only those terms which are positive.

\item{(ii)} $T$ acts on $\cW_d$ in an obvious way.
There is also a $\C^\times$ action on $\cW_d$ induced by
the standard action on $\P^1$ with weight $\alpha$. 
Each $\C^\times$
fixed point component in $\cW_d$ is (consisting of $K$-orbits of)
$$\eqalign{
\cY_r=\{\phi&=(x_1 w_0^{\bra D_1,r\ket} w_1^{\bra D_1,d-r\ket}
,...,x_N w_0^{\bra D_N,r\ket} w_1^{\bra D_N,d-r\ket})|\cr
&(x_1,..,x_N)\in\C^N,~
x_b=0~if~the~corresp.~monomial~has~negative~exponent\}.
}
$$
Let $j_r:\cY_r\ra \cW_d$ be the inclusion maps.
If nonempty, $\cY_r$ is canonically isomorphic to
a $T$-invariant submanifold in $X$ given by intersecting
those divisors $x_b=0$ corresponding to negative exponents above.
Denote the canonical inclusions by $\tau_r:\cY_r\ra X$.
Then ${\tau_r}_*(1)=\prod_{\bra D_a,r\ket<0~or~\bra D_a,d-r\ket<0}D_a$.
We will denote the class of
$D_a\cap \cY_r$ in $\cY_r$ simply by $D_a$.

\item{(iii)} The $G=\C^\times\times T$-equivariant 
Euler class of the normal bundle
of $\cY_r$ in $\cW_d$ is
$$
e_G(\cY_r/\cW_d)=
\prod_{\bra D_a,d\ket\geq0}
{\prod_{k=0}^{\bra D_a,d\ket}}_{k\neq\bra D_a,r\ket}
(D_a+\bra D_a,r\ket\alpha-k\alpha).
$$

\item{(iv)} Corresponding canonically to
every $T$-divisor class $D_a$ on $X$ is a 
$G$-divisor class $\hat D_a$ on $\cW_d$. It is determined by
$$j_r^*\hat D_a= D_a+\bra D_a,r\ket\alpha.$$
Similarly, every linear combination $D$ of the $D_a$ corresponds
to some $\hat D$ on $\cW_d$.

\lemma{$O(t)$ introduced above is a $1$-Euler series.}
\thmlab\ToricO
\proof
Let
\eqn\dumb{
\omega_d=
\prod_{\bra D_a,d\ket<0}
\prod_{k=1}^{-\bra D_a,d\ket-1}(\hat D_a+k\alpha)\in A^G_*(\cW_d).
}
By the localization theorem, 
$$
\int_{\cW_d} \omega_d\cdot e^{\hat H\cdot\zeta}
=\sum_r \int_{\cY_r}{j_r^*\omega_d\over e_G(\cY_r/\cW_d)}~
e^{(H+r\alpha)\cdot\zeta}.
$$
Obviously, the left hand side lies in $A_G^*(pt)[[\zeta]]\subset\cR[[\zeta]]$.
Now observe that the right hand side is
$$\eqalign{
\sum_r\int_{\cY_r}&
{\prod_{\bra D_a,d\ket<0}
\prod_{k=1}^{-\bra D_a,d\ket-1}
(D_a+\bra D_a,r\ket\alpha+k\alpha) \over
\prod_{\bra D_a,d\ket\geq0}
{\prod_{k=0}^{\bra D_a,d\ket}}_{k\neq\bra D_a,r\ket}
(D_a+\bra D_a,r\ket\alpha-k\alpha)}
~e^{(H+r\alpha)\cdot\zeta}\cr
&=\sum_r\int_X
{\prod_{\bra D_a,r\ket<0}
\prod_{k=0}^{-\bra D_a,r\ket-1}(D_a-k\alpha)\over
\prod_{\bra D_a,r\ket\geq0}
\prod_{k=1}^{\bra D_a,r\ket}(D_a+k\alpha)}\cr
&~~~\times {\prod_{\bra D_a,d-r\ket<0}
\prod_{k=0}^{-\bra D_a,d-r\ket-1}(D_a+k\alpha)\over
\prod_{\bra D_a,d-r\ket\geq0}
\prod_{k=1}^{\bra D_a,d-r\ket}(D_a-k\alpha)}.
~e^{(H+r\alpha)\cdot\zeta}\cr
&=\sum_r\int_X\bar O_r\cdot O_{d-r}
~e^{(H+r\alpha)\cdot\zeta}.
}
$$
This shows that $O(t)$ is a $1$-Euler series.
$\Box$

\remark{One can define the notion of Euler data on the basis of $\cW_d$
in a way similar to Definition \EulerData.
The classes \dumb~ in fact give an example of Euler data for $\cW_d$.
One can also construct the whole parallel theory of mirror principle 
for toric manifolds using $\cW_d$.}

\lemma{The two Euler series $O(t)$ and $\One(t)$ are linked.}
\thmlab\ToricLinking
\proof
Let $p\in X^T$, $pq$ be a balloon in $X$, $d=\delta[pq]\succ0$,
and $\lambda$ be the weight on the tangent line on $T_p(pq)$.
Let $F=(\P^1,f_\delta,0)$ be the fixed point in $M_{0,1}(d,X)$,
as given in Theorem \SimplePole, which says that
the function $i_p^*\One_d$ has the polar term,
at $\alpha=\lambda/\delta$,
\eqn\dumb{
e_T(p/X)~{1\over\lambda}{1\over\alpha-\lambda/\delta}
{1\over e_T(F/M_{0,1}(d,X))}.
}
We now compute the contribution from
$e_T(F/M_{0,1}(d,X))$ for a toric manifold $X$.
The virtual normal bundle of the point $F=(C=\P^1,f_\delta,0)$ in $M_{0,1}(d,X)$ is
$$
N_{F/M_{0,1}(d,X)}=[H^0(C,f_\delta^*TX)]-[H^1(C,f_\delta^*TX)]-A_C
$$
(notation as in section 2.3 \LLYI).
From the Euler sequence of $X$ \Jaczewski,
we get an equivariant exact sequence
for every balloon $pq$ in $X$,
$$0\ra \cO^{N-n}\ra \oplus_a \cO(D_a)|_{pq}\ra TX|_{pq}\ra 0$$
where $\cO$ is the trivial line bundle.
At $p$, there are exactly $n$ nonzero 
$D_a(p):=i_p^*D_a$ giving the weights
for the isotropic representation
$T_pX$, and $N-n$ zero $D_a(p)$ corresponding to
the trivial representation $\cO^{N-n}|_p$.
As usual we ignore the zero weights below, which 
must drop out at the end. 

Let $\lambda=D_b(p)$. Note that $\bra D_b,d\ket=1$ (section 2.3 \LLYII).
The bundle $\cO(D_b)$ contributes to 
$e_T([H^0(C,f_\delta^*TX)])$ the term
$$
\prod_{k=0}^{\delta-1}(D_b(p)-k\lambda/\delta).
$$
For each $a\neq b$ with $\bra D_a,d\ket\geq0$, 
the bundle $\cO(D_a)$ contributes to
$e_T([H^0(C,f_\delta^*TX)])$ the term
$$\eqalign{
&\prod_{k=0}^{\bra D_a,d\ket}(D_a(p)-k\lambda/\delta)~~if~D_a(p)\neq0,\cr
&\prod_{k=1}^{\bra D_a,d\ket}(D_a(p)-k\lambda/\delta)~~if~D_a(p)=0.
}$$
For each $a$ with $\bra D_a,d\ket<0$, 
the bundle $\cO(D_a)$ contributes to
$e_T([H^1(C,f_\delta^*TX)])$ the term
$$
\prod_{k=1}^{-\bra D_a,d\ket-1}(D_a(p)+k\lambda/\delta).
$$
The automorphism group $A_C$ contributes 
$e_T(A_C)=-\lambda/\delta$. 
Finally, we have
$$
e_T(p/X)=\prod_{D_a(p)\neq0}D_a(p).
$$
Combining all the contributions, we see that \dumb~ becomes
${1\over \alpha-\lambda/\delta}$ times
$$
{-1\over\delta}{\prod_{\bra D_a,d\ket<0}
\prod_{k=0}^{-\bra D_a,d\ket-1}(D_a(p)+k\lambda/\delta)\over
{\prod_{\bra D_a,d\ket\geq 0}}_{a\neq b}
\prod_{k=1}^{\bra D_a,d\ket}(D_a(p)-k\lambda/\delta)\times
\prod_{k=1}^{\delta-1}(D_b(p)-k\lambda/\delta)}.
$$
But this coincides with
$$lim_{\alpha\ra\lambda/\delta}~(\alpha-\lambda/\delta) i_p^*O_d.$$
This shows that 
$i_p^*O_d-i_p^*\One_d$ is regular at $\alpha=\lambda/\delta$.
$\Box$

Note that $O_d=\alpha^{\bra c_1(X),d\ket}+$ lower order terms, 
because $\sum D_a=c_1(X)$. Thus $O(t)$
is an Euler series linked to $\One(t)$ and 
meets the condition (**) of Theorem \ApplicationsII.
In particular to apply to the case $b_T=e_T$ or $c_T$,
all we need is the form \FactorizedForm~ for $\Omega$
and the bound \ChernClassBound. For then
Corollary \ApplicationsIII~ holds.

\bigskip
\noindent {\it Example.} Take $b_T=c_T$.
Take $V$ to be any direct sum
of convex equivariant line bundles $L_i$,
so that \ChernClassBound~ holds. Note that
in this case \FactorizedForm~ holds automatically. 
Then Theorem \ApplicationsII~ yields an explicit
formula for $A^{V,b_T}(t)$ in terms of the $O_d$ \FormulaO~
and the $P_d$ in Theorem \ChernPolynomial.
For $b_T=e_T$, and $V$ 
a direct sum of convex line bundles $L_i$ with $\sum_i c_1(L_i)=c_1(X)$,
we get a similar explicit formula for $A(t)$.
Plugging this formula into Theorem \KdTheorem
~ in the nonequivariant limit, we get

\corollary{Let
$$B(t)=e^{-H\cdot t/\alpha}\sum_d \prod_i\prod_{k=0}^{\bra c_1(L_i),d\ket}
(c_1(L_i)-k\alpha)\cap
{\prod_{\bra D_a,d\ket<0}
\prod_{k=0}^{-\bra D_a,d\ket-1}(D_a+k\alpha)\over
\prod_{\bra D_a,d\ket\geq0}
\prod_{k=1}^{\bra D_a,d\ket}(D_a-k\alpha)}~
e^{d\cdot t},
$$
as in Theorem \ApplicationsII.
Then we have
$$
\int_X\left(e^{f/\alpha}B(t)-e^{-H\cdot T/\alpha}\Omega\right)
=\alpha^{-3}(2\Phi-\sum T_i {\partial\Phi\over\partial T_i})
$$
where $T=t+g$, and $f,g$ are the
power series computed in Theorem \ApplicationsII.
}
This is
the general mirror formula in \HKTY\HLYII~ (see also references there),
formulated in the context of mirror symmetry
and reflexive polytopes \Batyrev\BB.

\subsec{A generalization}

We have now seen several ways to compute $A(t)=A^{V,b_T}(t)$
under various assumptions on either $V,b_T$, or $TX$,
or the classes $\One_d$, or some combinations of
these assumptions. We now combine these approaches to formulate an algorithm
for computing $A(t)$ in full generality on any balloon manifold $X$,
for arbitrary $V,b_T$. The result will be in terms of
certain (computable) $T$ representations.

\item{(i)} By Lemma \SquareLemmaIII, the $A_d$ is of the form
$$A_d={\tau^*\phi_d\over e_G(X/W_d)},$$
where $\phi_d\in A^T_*(Y)[\alpha]$, 
hence represented by a polynomial $\C[\cT^*][H_1,...,H_m,\alpha]$.
Note that the denominator of $A_d$ is
$e_G(X/W_d)=e_T(X/Y)\cdot \tau^* e_G(Y_0/W_d)$.
Thus the goal is to compute the class $\tau^*\phi_d$ for all $d$.
We shall set up a (over-determined) system of
linear equations with a solution (unique up to $ker~\tau^*$)
given by the $\phi_d$.

\item{(ii)} By Theorem \KdTheorem, the degree of 
$\phi_d$ is bounded according to 
$$deg_\alpha A_d\leq -2.$$

\item{(iii)} By Theorem \SimplePole~(i), at any fixed point $p$,
the function $i_p^*A_d$ is
regular away from $\alpha=0$ or $\lambda/\delta$,
where $\lambda$ is a weight on $T_pX$. In other words,
$$Res_{\alpha=\gamma}(\alpha-\gamma)^k~i_p^* A_d=0$$
for all $\gamma\neq 0,\lambda/\delta$ and $k\geq0$.
Note that these are all linear conditions on $\phi_d$.

\item{(iv)} By Theorem \SimplePole~ (ii) (see notations there), 
for any balloon $pq$ in $X$ and $d=\delta[pq]\succ0$, we have
\eqn\dumb{\eqalign{
lim_{\alpha\ra\lambda/\delta}(\alpha-\lambda/\delta)~i_p^*A_d
&={e_T(p/X)\over\lambda~ e_T(F/M_{0,1}(d,X))}
i_F^*\rho^*b_T(V_d)\cr
&={-e_T(p/X)\over\delta}{e_T[H^1(\P^1,f^*_\delta TX)]'
\over e_T[H^0(\P^1,f^*_\delta TX)]'}~
b_T(i_{\rho(F)}^*V_d).
}}
Here we have used the fact that
$N_{F/M_{0,1}(d,X)}=[H^0(\P^1,f_\delta^*TX)]-[H^1(\P^1,f_\delta^*TX)]-A_C$
(cf. section 5.2). The prime in the Euler classes above
means that we drop the zero weights in the
$T$ representations
$[H^i(\P^1,f^*_\delta TX)]$.
Now if $V=V^+\oplus V^-$ is a convex/concave bundle on $X$, then
we have the $T$ representation
$$
i_{\rho(F)}^*V_d=H^0(\P^1,f_\delta^*V^+)\oplus H^1(\P^1,f_\delta^* V^-).
$$
Thus
$$
b_T(i_{\rho(F)}^*V_d))=
b_T[H^0(\P^1,f_\delta^*V^+)]~b_T[H^1(\P^1,f_\delta^* V^-)],
$$
which is just the value of $b_T$ for a trivial bundle over a point.
Note that if $U$ is any $T$ representation with weight decomposition
$U=\oplus_i \C_{\nu_i}$, then 
$b_T(U)=\prod_i b_T(\C_{\nu_i})$
by the multiplicativity of $b_T$.
Hence once the $T$ representations $H^i(\P^1,f_\delta^*V^\pm)$
are given,  eqn. \dumb~
becomes a linear condition on the $\phi_d$, where the right hand
side is some known element in $\C(\cT^*)$.

\item{(v)} Finally, we know that $A(t)$ is an $\Omega$-Euler series.
This is (inductively) a linear condition on the $\phi_d$.

\item{(vi)} By Theorem \Uniqueness, any
solution to the linear conditions in (ii)-(v) necessarily represents
the class $\tau^*\phi_d$ we seek.

\bs
Of course, this algorithm relies on knowing
the $T$ representations  $[H^i(\P^1,f_\delta^*TX)],~
[H^i(\P^1,f_\delta^*V)]$ induced by the
$T$-equivariant  bundles $TX|_{pq}$ and $V|_{pq}$
on each balloon $pq\cong\P^1$.
But describing them for any given $X$ and $V$ is
clearly a classical question.
We have seen that these representations are
easily computable in many cases.
We now discuss a general situation in which 
these representations can also be computed 
similarly.

Let $V$ be any $T$-equivariant vector bundle on $X$ and let
$$0\ra V_N\ra\cdots\ra V_1\ra V\ra 0$$
be an equivariant resolution. Then by the Euler-Poincare
Principle, we have
$$[H^0(\P^1,f_\delta^* V)]-[H^1(\P^1,f_\delta^* V)]
=\sum_a(-1)^{a+1} 
([H^0(\P^1,f_\delta^* V_a)]-[H^1(\P^1,f_\delta^* V_a)]).
$$
Note that
there is a similar equality of 
representations whenever $V$ is a term
in any given exact sequence
$$0\ra V_N\ra\cdots\ra V_i\ra V\ra V_{i-1}\ra\cdots\ra V_1\ra 0.$$
Now suppose that each $V_a$ is a
direct sum of $T$-equivariant line bundles.
Then each summand $L$ will contribute to
$[H^0(\P^1,f_\delta^* V_a)]-[H^1(\P^1,f_\delta^* V_a)]$
the representations 
$$c_1(L)|_p-k\lambda/\delta,~~~k=0,1,...,l\delta,$$
if $l=\bra c_1(L),[pq]\ket\geq0$; and
$$c_1(L)|_p+k\lambda/\delta,~~~k=1,...,l\delta-1,$$
if $-l=\bra c_1(L),[pq]\ket<0$ (cf. proof of Theorem \EulerClass).
In this case, 
$[H^0(\P^1,f_\delta^* V)]-[H^1(\P^1,f_\delta^* V)]$ are
then determined completely.
Thus whenever a $T$-equivariant resolution
by line bundles is known for $TX$ and the convex/concave bundle $V^\pm$,
the right hand side of eqn. \dumb~ becomes known.

\bigskip
\noindent{\it Example.} Consider the case
$X=\P^n$, $V=TX$, and $b_T$ the Chern polynomial.
This will be an example where $V$ has no splitting type,
but $A(t)$ can be computed
via a $T$-equivariant resolution nevertheless.
Recall the $T$ equivariant Euler sequence
$$
0\ra\cO\ra\oplus_{i=0}^n\cO(H-\lambda_i)\ra TX\ra0.
$$
For $F=(\P^1,f_\delta,0)$, where $f_\delta$ is
the $\delta$-fold cover of the balloon $pq$, this gives
$$
b_T(i^*_{\rho(F)}V_d)={1\over x}
\prod_i\prod_{k=0}^d(x+\lambda_j-\lambda_i-k\lambda/\delta).
$$
Here $p,q$ are the $j$th and the $l$th fixed points in
$\P^n$, so that $\lambda=\lambda_j-\lambda_l$.
We can use this to set up a linear system to solve for $A(t)$
inductively. However, there is an easier way to compute $A(t)$
in this case.
Observe that
$\Omega=b_T(V)={1\over x}\prod_i(x+H-\lambda_i)$, and that
$$P:~~~P_d:={1\over x}\prod_i\prod_{k=0}^d(x+\kappa-\lambda_i-k\alpha)$$
defines an $\Omega$-Euler data (see section 2.2 \LLYI). 
Since $j_0^*\kappa=H$ and $i_p^*H=\lambda_j$, it
follows that
$$b_T(i^*_{\rho(F)}V_d)=i_p^*j_0^*P_d$$
at $\alpha=\lambda/\delta$.
By the corollary to Theorem \SimplePole, the Euler series
$$B(t):=e^{-H\cdot t/\alpha}\sum j_0^*P_d\cap\One_d$$
is linked to $A(t)$. Obviously, we have $deg_\alpha j_0^*P_d=(n+1)d$,
hence $P$ meets condition (*) of Theorem \ApplicationsII
~($\tau$ is the identity map here).
For $O_d=\One_d$, condition (**) there is also automatic.
It follows that
$$A(t+g)=e^{f/\alpha} B(t)$$
where $f,g$ are explicitly computable functions
from Theorem \ApplicationsII.
Note that $rank~V_d=(n+1)d+n$, and so Theorem \KdTheoremII~
yields immediately the codimension $s=3$ Chern class of $V_d$.

\subsec{Blowing up the image}

In this section, we discuss another approach to compute $A(t)$.
For clarity, we restrict to the case
of a convex $T$-manifold $X$ ($T$ may be trivial), 
and $b_T\equiv1$, so that
$A(t)=\One(t)$. 
Thus we will study the classes
$$
\One_d=e^X_*\left({LT_{0,1}(d,X)\over e_G(F_0/M_d(X))}\right).
$$
We will actually be interested in the 
integrals $\int_X \tau^*e^{H\cdot\zeta}\cap\One_d$,
where $\tau:X\ra Y$ is a given projective embedding.
For the purpose of studying the intersection numbers in section 3.3,
this is adequate.
Since $X$ is assumed convex, $LT_{0,1}(d,X)$ is
represented by $M_{0,1}(d,X)$. Likewise for $M_d(X)$.

Suppose that 
we have a commutative diagram
\eqn\CommSquareIII{
\matrix{
F_0 & {\br e^Y\over \longrightarrow} & Y_0
& {\br g\over\longleftarrow} & E_0\cr
\downarrow i &  & \downarrow j
&   &\downarrow k\cr
M_d & {\br \varphi\over \longrightarrow} & W_d
& {\br \psi\over \longleftarrow} & \cQ_d.}
}
Here the left hand square is as in \CommSquareII~
($i_0,j_0,M_d(X)$ there are written as $i,j,M_d$ here for clarity).
We assume that 
$\cQ_d$ is a $G$-manifold,
that $\psi:\cQ_d\ra W_d$ is a $G$-equivariant resolution of
singularities of $\varphi(M_d)$, and that $E_0$ is
the fixed points in $\psi^{-1}(Y_0)$.
Here $g$ denotes the restriction of $\psi$,
and $k$ the inclusion.
Recall that $\varphi$ is an isomorphism into its
image away from the singular locus of $M_d$.
The singularities in $\varphi(M_d)$ is
the image of the compactifying divisor in $M_d$,
which has codimension at least $2$.
Then we have the equality in $A_*^G(W_d)$:
$$
\varphi_*[M_d]=\psi_*[\cQ_d].
$$

Applying functorial localization to the left hand square
in \CommSquareIII~ as in section 3.1, we get
$$
{j^*\varphi_*[M_d]\over e_G(Y_0/W_d)}=
e^Y_*\left({[F_0]\over e_G(F_0/M_d)}\right).
$$
Doing the same for the right hand square, we get
$$
{j^*\psi_*[\cQ_d]\over e_G(Y_0/W_d)}
=g_*\left({[E_0]\over e_G(E_0/\cQ_d)}\right).
$$
It follows that

\lemma{ 
In $A_*^T(Y)$, we have the equality
$$
\tau_*\One_d=e^Y_*\left({[F_0]\over e_G(F_0/M_d)}\right)
=g_*\left({[E_0]\over e_G(E_0/\cQ_d)}\right).
$$
}
It follows that
$$\eqalign{
\int_X \tau^*e^{H\cdot \zeta}\cap\One_d
&=\int_{Y_0} e^{H\cdot\zeta}\cap
g_*\left({[E_0]\over e_G(E_0/\cQ_d)}\right)\cr
&=\int_{E_0} {g^*e^{H\cdot\zeta}
\over e_G(E_0/\cQ_d)}\cr
&=\int_{E_0} {g^*j^*e^{\kappa\cdot\zeta}
\over e_G(E_0/\cQ_d)}\cr
&=\int_{E_0} {k^*\psi^*e^{\kappa\cdot\zeta}
\over e_G(E_0/\cQ_d)}.
}$$

In many cases, the spaces $\cQ_d$ can be explicitly described, and 
the classes $\psi^*\kappa$ on $\cQ_d$ can be expressed in terms of  
certain universal classes. 
For example, when $X$ is a flag variety, 
then the $\cQ_d$ 
can be chosen to be the Grothendieck Quot scheme (cf. \CF\Be). 
Integration on the Quot scheme can be done by explicit
localization (cf. \Stromme). When $X$ is a Grassmannian and $\tau:X\ra Y=\P^N$
is the Plucker embedding, then
$\psi^*\kappa=-c_1(S)$, where $S$ is a universal subbundle 
on the Quot scheme. In this case, the image $\psi(Q_d)$
has been studied extensively in \So\SS. 

When $X$ is not convex, 
a similar method still works if we can find 
an explicit cycles $Z_d$ in $\cQ$ such that
$$
\varphi_*LT_d(X)=\psi_*[Z_d]
$$
in $A_*^G(W_d)$.
This approach deserves further investigation.

\subsec{Higher genus}

In this section, we discuss a generalization of mirror principle
to higher genus. More details will appear elsewhere.
As before $X$ will be a projective $T$-manifold,
and $\tau:X\ra Y$ a given $T$-equivariant projective embedding.
(Again, $T$ may be the trivial group.)

Let $M_{g,k}(d,X)$ denote the $k$-pointed, arithmetic genus $g$,
degree $d$, stable map moduli stack of $X$. Let $M_d^g$ denote 
$M_{g,0}((d,1),X\times\P^1)$. Note that
for each stable map $(C,f)\in M_d^g$
there is a unique branch $C_0\cong\P^1$ in $C$
such that $f$ composed with the projection $X\times\P^1\ra \P^1$
maps $C_0\ra\P^1$ isomorphically. Moreover, $C$ 
is a union of $C_0$ with some disjoint curves $C_1,..,C_N$, 
where each $C_i$ intersects $C_0$ at a point $x_i\in C_0$.
The map $f$ composed with
$X\times\P^1\ra \P^1$ collapses all $C_1,..,C_N$.

The standard $\C^\times$ action on $\P^1$
induces an action on $M_d^g$. The fixed point components are
labelled by 
$F_{d_1,d_2}^{g_1,g_2}$ with $d_1+d_2=d$, $g_1+g_2=g$. 
As in the genus zero case, 
a stable map $(C,f)$ in this component is given by gluing
two $1$-pointed stable maps $(f_1,C_1,x_1)\in M_{g_1,1}(d_1,X),~ 
(f_2,C_2,x_2)\in M_{g_2,1}(d_2,X)$ with $f_1(x_1)=f_2(x_2)$,
to a $\P^1$ at $0$ and $\infty$ at the marked points (cf. section 3).
We can therefore identify 
$F_{d_1,d_2}^{g_1,g_2}$ with
$M_{g_1,1}(d_1,X)\times_X M_{g_2,1}(d_2,X)$.
We denote by
$$F_{d_1,d_2}^g:=
\coprod_{g_1+g_2=g} F_{d_1,d_2}^{g_1,g_2},~~~
i_{d_1,d_2}:F_{d_1,d_2}^g\ra M_d^g,
$$
the disjoint union and inclusions.
There are two obvious projection maps
$$
p_0:F^g_{d_1,d_2}\ra \coprod_{g_1=0}^g M_{g_1,1}(d_1,X),~~
p_\infty:F^g_{d_1,d_2}\ra \coprod_{g_2=0}^g M_{g_2,1}(d_2,X).
$$
The map $p_0$ strips away the stable maps $(f_2,C_2,x_2)$
glued to $\infty$ and forgets the $\P^1$;
$p_\infty$ strips away the stable map
$(f_1,C_1,x_1)$ glued to $0$ and forgets the $\P^1$.
We also have the usual evaluation maps, and the forgetting map:
$$e_{d_1,d_2}:F^g_{d_1,d_2}\ra X,~~~e_d:M_{g,1}(d,X)\ra X,~~~
\rho:M_{g,1}(d,X)\ra M_{g,0}(d,X).$$
Relating and summarizing the natural maps above is the following
diagram:
\eqn\CommDiagram{
\matrix{
X & {\br e_{d_1,d_2}\over\lla} &  F^{g_1,g_2}_{d_1,d_2} &
{\br i_{d_1,d_2}\over\lra} &  M_d^g & {\br \pi\over\lra} & M_{g,0}(d,X) \cr
e_{d_1}\uparrow & p_0\swarrow &   &\searrow p_\infty  & & & \cr
M_{g_1,1}(d_1,X) &   &   &    & M_{g_2,1}(d_2,X)  &  & \cr 
\rho\downarrow &  &   &  &\downarrow\rho &  &  \cr
M_{g_1,0}(d_1,X) &   &   &    & M_{g_2,0}(d_2,X)  &  & 
}
}

Fix a class $\Omega\in A_T^*(X)$.
We call a list $b_d^g\in A_T^*(M_{g,0}(d,X))$ an $\Omega$-gluing sequence
if we have the identities on the $F^g_{d_1,d_2}$:
$$
e_{d_1,d_2}^*\Omega\cdot i_{d_1,d_2}^*\pi^* b_d^g
=\sum_{g_1+g_2=g} p^*_0\rho^* b_{d_1}^{g_1}\cdot
p^*_\infty\rho^* b_{d_2}^{g_2}.
$$
It is easy to verify that $b_d^g\equiv1$ is an example of a $1$-gluing
sequence.
Restricted to $g=0$, the identity above is precisely the gluing identity
in section 3.2. There we have found that the gluing identity
results in an Euler series. It turns out that
a gluing sequence too leads to an Euler series.
For $\omega\in A_G^*(M_d^g)$ and $d=d_1+d_2$, define (cf. section 3.2)
$$
i^{vir}_{d_1,d_2}\omega:={e_{d_1,d_2}}_*
\left({i_{d_1,d_2}^*\omega\cap [F^g_{d_1,d_2}]^{vir}
\over e_G(F^g_{d_1,d_2}/M_d^g)}\right)\in A_*^T(X)(\alpha).
$$
Then for a given gluing sequence $b_d^g\in A_T^*(M_{g,0}(d,X))$, we have the
identities
$$
\Omega\cap 
i^{vir}_{d_1,d_2}\pi^*b_d^g
=\sum_{g_1+g_2=g}
\overline{i_{0,d_1}^{vir}\pi^*b_{d_1}^{g_1} } 
\cdot i_{0,d_2}^{vir}\pi^*b_{d_2}^{g_2}. 
$$
Again, putting $g=0$, we get the identities in Theorem \GluingTheorem.
The argument in the higher genus case is essentially
the same as the genus zero case. Here, one chases through
a fiber diagram analogous to \FiberSquare~ using the
associated refined Gysin homomorphism, together with
the diagram \CommDiagram.

Now given a gluing sequence, we put
$$
A_d^g:=i_{0,d}^{vir}\pi^*b_d^g,~~~ A_d:=\sum_g A_d^g~ \mu^g,~~~
A(t):=e^{-H\cdot t/\alpha}\sum_d A_d~e^{d\cdot t}.
$$
Here $\mu$ is a formal variable.
Then $A(t)$ is an Euler series.
(We must, of course, replace the ring $\cR$ by
$\cR[[\mu]]$.)
The argument is also similar to
the genus zero case: one applies
functorial localization to the diagram
$$\matrix{
F^g_{d_1,d_2} & {\br i_{d_1,d_2}\over \lra} & M_d^g\cr
e_{d_1,d_2}\downarrow &  & \downarrow \varphi\cr
~~X\subset Y & {\br j_{d_1,d_2}\over \lra} & W_d},
$$
the same way we have done to diagram \CommSquareII~ in section 3.2.

We can proceed further in a way parallel
to the genus zero case. Namely, to find further constraints
to a gluing sequence, we should compute the linking
values of the Euler series $A(t)$. For this, let's assume
that $X$ is a balloon manifold, as in sections 4.1 and 5.3. 
In genus zero, the linking values of an Euler series,
say coming from $b_T(V_d)$, are determined by the
restrictions $i_F^*b_T(\rho^* V_d)$ to the
isolated fixed point $F=(\P^1,f_\delta,0)\in M_{0,1}(d,X)$,
which is a $\delta$-fold cover of a balloon $pq$ in $X$
(see Theorem \SimplePole).
In higher genus, this is replaced by
a component in $M_{g,1}(d,X)$ consisting of
the following stable maps $(C,f,x)$.
Here $C$ is a union of two curves $C_1$ and
$C_0\cong\P^1$ such that $C_0{\br f\over\ra} pq$ is a $\delta$-fold
cover with $f(x)=p$, and $f(C_1)=q$. 
Therefore this fixed point component can be 
identified canonically with
$\bar M_{g,1}$, the moduli space of genus $g$, $1$-pointed, stable curves.
Let's abbreviate it $F$.
The linking values of $A(t)$ for this component
is then a power series summing over integrals on
$\bar M_{g,1}$ of classes  given in terms of $i_F^*\rho^* b_d^g$
and $e_T(F/M_{g,1}(d,X))$ (cf. Theorem \SimplePole).

\footatend\vfill\supereject\immediate\closeout\rfile\writestoppt
\baselineskip=14pt\centerline{{\bf References}}\bigskip{\frenchspacing%
\parindent=20pt\escapechar=` \input refs.tmp\vfill\eject}\nonfrenchspacing

\end